%% file: An_Equivariant_Framed_Pontryagin-Thom_Theorem_for_Families_of_Manifolds.tex
\documentclass[11pt]{amsart}
\title{The Pontryagin-Thom Theorem for Families of Framed Equivariant Manifolds}
\author{Lucas Williams}
\address{Department of Mathematics and Statistics, Binghamton University, Binghamton, NY 13902, USA}
\email{lwilli39@binghamton.edu}

\usepackage{amsmath,amssymb,amsthm,array,color,multirow,pdflscape,mathrsfs,subcaption,todonotes}
\usepackage{tikz,tikz-cd}
\usetikzlibrary{bending}
\usepackage{float}
\usepackage{imakeidx}
\usepackage[toc,page]{appendix}

\usepackage[margin=1.25in]{geometry}
\newcommand{\beforesubsection}{\vspace{1em}}
\newcommand{\aftersubsection}{}

\include{preamble_defs}
\include{preamble_refs_angelica}
\begin{document}

\maketitle

\begin{abstract}
	The Pontryagin-Thom theorem gives an isomorphism from the cobordism group of framed $n$-manifolds to the $n$th stable homotopy group of the sphere spectrum. In this paper, we prove the generalization of the Pontryagin-Thom theorem for families of framed equivariant manifolds parameterized over a compact base space.
\end{abstract}

\setcounter{tocdepth}{2}
\tableofcontents

\parskip 2ex

\section{Introduction}

A framed manifold is a smooth compact manifold along with a trivialization of its stable tangent bundle. The set of framed $n$-manifolds up to cobordism forms an abelian group under disjoint union which we denote $\omega_n$. The Pontryagin-Thom theorem, also known as the Pontryagin-Thom isomorphism, gives a correspondence between cobordism classes of framed manifolds and the stable homotopy groups of spheres:
\[
\omega_n \cong\pi_n(\Sph).
\]
This theorem provides a bridge between geometric topology and homotopy theory. The Pontryagin-Thom theorem may be seen as the jumping off point for a body of work that has led to such celebrated results as the solution to the Kervaire invariant one problem \cite{hhr}, various computations of smooth structures on spheres as in \cite{kervaire-milnor} and \cite{wang_xu_61stem}, and the surgery theoretic classification of high dimensional manifolds \cite{luck_surgery}. In the present work, we provide a new and more general proof of the equivariant Pontryagin-Thom theorem. 

\begin{thm}\label{thm:first main theorem}
Let $G$ be a product of a finite group and a torus. There is an isomorphism 
\[
\omega_V^G(X,A)\to \pi_V^G(\Sigma^\infty X/A)
\]
from the group of $V$-framed $G$-manifolds admitting a continuous equivariant map of pairs $(M,\partial M)\to (X,A)$ to the $V$th equivariant stable homotopy group of the suspension spectrum of $X/A$. Here a $V$-framed manifold is a manifold $M$ equipped with an equivalence class of stable bundle isomorphism between $TM$ and $M\times V$. The precise definition will be given in \cref{df:V-framing}.
\end{thm} 

Our treatment of this theorem differs from the classical literature in several respects. We use the modern framework of genuine $G$-spectra and the language of $G$-manifolds with corners. We also generalize the group of equivariance from a finite group (as is the case in \cite{kosniowski_bordism}) to any product of a finite group and a torus. As noted in \cite[6.2.13]{schwede_global} this is the largest class of groups for which such a result should be expected to hold.

The main advantage of our new proof of the equivariant Pontryagin-Thom theorem for framed $G$-manifolds is that it generalizes quite nicely to the fiberwise setting. This generalization is our second main theorem. 

\begin{thm}\label{thm:second main theorem}
Let $G$ be a product of a finite group and a torus. There is an isomorphism 
\[
\omega_V^G((X,A)\to B) \to \pi_V^G(\Gamma_BR^{lv}\Sigma^\infty X\cup_A B)
\]
from the cobordism group of families of $V$-framed $G$-manifolds parameterized over an equivariant map $(X,A)\to B$ and the $V$th equivariant stable homotopy group of the right derived spectrum of sections of the suspension spectrum of $X\cup_A B$. 
\end{thm}

In the non-equivariant setting, given an element $(S^{k+n}\xrightarrow{f} S^k)\in \pi_n(\Sph)$, a framed manifold is produced by deforming $f$ to be smooth and transverse to $0\in S^k$, then taking $f^{-1}(0)$. In the equivariant setting it is not the case that a map between manifolds may be deformed so as to be transvere to a submanifold of our choice. Due to this lack of a `transverse approximation theorem' it is frequently (erroneously) claimed that the Pontryagin-Thom theorem fails for framed $G$-manifolds. 

However, results from as early as the 1970's prove a Pontryagin-Thom style theorem for $G$-manifolds in both the framed and unoriented cases. Segal originally announced a Pontryagin-Thom theorem for framed $G$-manifolds in \cite{segal1970equivariant}. The proof was completed by Hausschild \cite{Hausschild1974} and Kosniowski \cite{kosniowski_bordism}. Pontryagin-Thom theorems for (un)oriented $G$-manifolds appear in tom Dieck \cite{tomDieckEquivariantBordism}, Wasserman \cite{wasserman}, and Waner \cite{waner_1984_eq_bordism}. See Schwede \cite{schwede_global} for a modern treatment of tom Dieck's result.

\subsection*{Organization}
In section two we develop the theory of cobordism for framed $G$-manifolds before proving \cref{thm:first main theorem}. 

In section three we develop the theory of cobordism for families of framed $G$-manifolds and then prove \cref{thm:second main theorem}. The proof relies largely on equivariant parameterized stable homotopy theory.

\textbf{Acknowledgements}
The author would like to acknowledge John Klein, Cary Malkiewich, and Kate Ponto for many helpful conversations. The author would also like to thank several anonymous reviewers whose suggestions improved this paper. The author was partially supported by NSF grant DMS-2052923. This paper represents a portion of the author's thesis written under the direction of Cary Malkiewich at Binghamton University.

\section{The equivariant version}\label{sec:Equivariant PT}

In this section we prove a Pontryagin-Thom theorem for framed $G$-manifolds. Many of the ideas in this section originate in \cite{kosniowski_bordism}. However, unlike \cite{kosniowski_bordism} we use the modern framework of genuine $G$-spectra as well as the theory of $G$-manifolds with corners. Moreover, while the group of equivariance in \cite{kosniowski_bordism} is always finite, we will generalize the result to the product of a finite group with a torus. 

\beforesubsection
\subsection{Preliminaries}\aftersubsection

\beforesubsection
\subsubsection{Equivariant framed bordism}\label{subsec:equivariant cobordism}\aftersubsection

Throughout this paper, $G$ will be a product of a finite group with a torus and a $G$-representation will mean a finite dimensional orthogonal real representation of $G$. All manifolds are compact (possibly with boundary).

\begin{notn}
If $V$ is a $G$-representation, and $M$ is a $G$-manifold then write $\varepsilon_M(V)$ to denote the product bundle $M\times V$.  
\end{notn}

\begin{df}
A \textbf{$G$-manifold} is a compact smooth manifold equipped with a smooth $G$-action. 
\end{df}

\begin{rmk}
We will use the term `closed manifold' to mean a compact smooth manifold without boundary. When we refer to a compact manifold rather than a closed manifold, we are allowing the possibility of a boundary. 
\end{rmk}

\begin{df}\label{df:V-framing}
Let $M$ be a $G$-manifold and $V$ a $G$-representation. A $V$-framing on $M$ is an equivalence class of $G$-vector bundle isomorphisms
\begin{equation}\label{eq:V-framing}
TM\oplus \varepsilon_M(\R^n) \cong \varepsilon_M(V)\oplus\varepsilon_M(\R^n)
\end{equation}
We say $G$-vector bundle isomorphisms $\phi$ and $\psi$ as in \cref{eq:V-framing} are equivalent if they are $G$-homotopic over $M$. 

We also say that
\[
TM\oplus \varepsilon_M(\R^n) \xrightarrow{\phi} \varepsilon_M(V)\oplus\varepsilon_M(\R^n)
\]
is equivalent to 
\[
TM\oplus \varepsilon_M(\R^{n+1}) \xrightarrow{\psi} \varepsilon_M(V)\oplus\varepsilon_M(\R^{n+1})
\]
if $\psi$ is obtained from $\phi$ by extending $\phi$ to be the identity in the $(n+1)$st coordinate of each fiber of $\varepsilon_M(\R^{n+1})$. 
\end{df}

Note that although the manifolds in the above definition are smooth, we are working with continuous vector bundles over them.

\begin{rmk}
There is another notion of a $V$-framed $G$-manifold, in which we replace $\R^n$ in the above definition with any finite dimensional real orthogonal $G$-representation. While these two notions are equivalent if $M$ is a free $G$-manifold, our definition of a framing is, in general, stricter \cite{kosniowski_bordism}. In \cite{waner_1984_eq_bordism}, the author considers bordism of framed $G$-manifolds using this second notion of framing. This results in the same representing spectrum as in our notion of framed cobordism. Thus, while the opposing definitions of a framing are certainly not equivalent, they result in isomorphic cobordism theories. 
\end{rmk}

\begin{df}\label{df:pair}
A \textbf{$G$-topological pair} $(X,A)$ is a $G$-space $X$ together with a $G$-invariant subspace $A$.
\end{df}

\begin{df}
Let $(X,A)$ be a $G$-topological pair and $V$ a $G$-representation. A \textbf{$V$-framed bordism element of $(X,A)$} is a pair $(M,f)$ such that
\begin{enumerate}
\item $M$ is a $V$-framed $G$-manifold. 
\item $f:(M,\partial M)\to (X,A)$ is a continuous equivariant map of pairs. In other words, $f:M\to X$ is such that $f(\partial M)\subset A$. 
\end{enumerate}
\end{df}

Before defining the cobordism relation, we will need a few more definitions.

If $M$ is a $V$-framed $G$-manifold then it comes equipped with an isomorphism
\[
TM\oplus\varepsilon_M(\R^n)\xrightarrow{\phi} \varepsilon_M(V)\oplus \varepsilon_M(\R^n).
\]
Let $-\phi$ denote the trivialization
\[
-\phi = \phi\oplus -1_\R: TM\oplus \varepsilon_M(\R^n)\oplus\varepsilon_M(\R) \cong \varepsilon_M(V)\oplus \varepsilon_M(\R^n)\oplus \varepsilon_M(\R)
\]
where $-1_\R:M\times \R\to M\times\R$ is given by $(m,t)\mapsto (m,-t)$. 

Call the manifold with this new framing $-M$. 

Before defining equivariant framed cobordisms, we must first define an appropriate notion of $G$-manifolds with corners. The next several definitions build up to this.

\begin{df}
An $n$-dimensional manifold is a \textbf{smooth manifold with depth-two corners} if it is a second countable Hausdorff space with a maximal smooth atlas locally modeled on $[0,\infty)^2\times \R^{n-2}$. We will shorten this to \textbf{smooth manifold with corners} from here forward. If $M$ is a smooth manifold with corners, then $x\in M$ is a \textbf{corner point} if there is a chart identifying $x$ with the origin of $[0,\infty)^2\times\R^{n-2}$. The set of corner points of $M$ is denoted $\partial^cM$. 
\end{df}

\begin{df}
If $M$ is a smooth manifold with corners and $x\in M$, the set of \textbf{local boundary components near $x$} is the following inverse limit over neighborhoods $U$ of $x$:
\[
\underset{x\in U}{\lim}\pi_0(U\cap \partial M\backslash \partial^c M).
\]
\end{df}

\begin{df}
Let $M$ be a smooth manifold with corners and $G$ a compact lie group acting on $M$. The action of $G$ is \textbf{trivial along corners} if either of the following equivalent conditions hold. 
\begin{enumerate}
\item For each point $x\in \partial M$, the action of the stabilizer $G_x$ on the set of local boundary components near $x$ is trivial. 
\item $M$ is locally modeled on $G\times_H\times V\times [0,\infty)^2$ for varying $H\leq G$ and $H$-representations $V$. 
\end{enumerate}
\end{df}
For a proof of these equivalences and for (non)examples of groups acting trivially on corners, see \cite[2.6]{gimm}.

\begin{df}\label{df:G-manifold with corners}
A \textbf{smooth $G$-manifold with corners} is a smooth manifold with corners equipped with a group action that is trivial along corners.
\end{df}

\begin{df}
A function on an open subset of $\R^{n-2}\times[0,\infty)^2$ is smooth if it extends to a smooth function on an open subset of $\R^n$. A function on a smooth manifold with corners is smooth if it is smooth in the manifolds coordinate charts. 
\end{df}

This definition allows us to define the tangent bundle of a smooth manifold with corners in the usual way using derivations of smooth functions.

\begin{df}
Let $M$ be a smooth manifold with corners. The \textbf{smooth boundary of $M$}, denoted $\tilde{\partial}M$, is the set of all pairs $(x,b)$ where $x\in M$ and $b$ is a local boundary component near $x$. The smooth boundary inherits a smooth structure from $M$. The map $i:\tilde{\partial}M\to M$ is smooth and its image is $\partial M$. It fails to be injective precisely at the corners of $M$.
\end{df}

\begin{df}\label{df:face}
Let $M$ be a $G$-manifold with corners. A \textbf{face} $F$ of $M$ is a $G$-invariant subspace of the smooth boundary $\tilde{\partial}M$ such that $F$ is a union of components of $\tilde{\partial} M$ and the map $i:\tilde{\partial}M\to \partial M$ is injective when restricted to $F$. 
\end{df}

\begin{df}\label{df:G-cobordism}
Let $M$ and $N$ be smooth $n$-dimensional $G$-manifolds with boundary. A \textbf{$G$-cobordism between $M$ and $N$} is an $(n+1)$-dimensional $G$-manifold with corners, $W$, equipped with the following structure. There is a face of $W$ identified with $M$ called the bottom and a face of $W$ identified with $N$ called the top. The closure of the complement of $M\amalg N$ in $\partial W$ is called the sides. We also require that the corners of $W$ coincide exactly with $\partial M\amalg \partial N$.
\end{df}

\begin{rmk}
Given a cobordism $W$ from $M$ to $N$, the sides of $W$ give a cobordism from $\partial M$ to $\partial N$. 
\end{rmk}

\begin{df}\label{df:cobordism}
We will say that two $V$-framed cobordism elements $(M,f)$ and $(N,g)$ of $(X,A)$ are \textbf{$V$-framed cobordant} if there exists some pair $(W,q)$ such that 
\begin{description}
\item[(i)] $W$ is a $(V\oplus \R)$-framed $G$-manifold with corners.
\item[(ii)] $W$ is a $G$-cobordism between $M$ and $N$ as in \cref{df:G-cobordism}. 
\item[(iii)] $M\amalg N\subset \partial W$ such that the induced $V$-framing on $\partial W$ restricted to $M \amalg N$ agrees with the framing on $M\amalg N$.
\item[(iv)] $q:W\to X$ is an equivariant map such that $q|_{M}=f$ and $q|_{N}=g$ and $q$ maps the sides of $W$ into $A$. 
\end{description}
We will frequently shorten $V$-framed cobordism to cobordism when the context is clear.
\end{df}

Regarding point (iii) in \cref{df:cobordism}, recall from \cref{df:G-cobordism} that we identity $M$ and $N$ with faces of $W$ called the bottom (denoted $F_0$) and top (denoted $F_1$) respectively. Recall from \cref{df:face} that $F_0$ and $F_1$ both inject into $W$. Also recall that a framing on $W$ is an equivalence class of $G$-bundle isomorphisms
\[
TW\oplus \varepsilon_W(\R^k) \cong \varepsilon_W(V\oplus\R)\oplus\varepsilon_W(\R^k)
\]
We induce a framing on $F_0$ from the framing on $W$ by pulling back the framing on $W$ along the inclusion $F_0\hookrightarrow W$. We induce a framing on $F_1$ from the framing on $W$ by pulling back the framing on $W$ along $F_1\hookrightarrow W$ and then negating the $\R$ coordinate of $\varepsilon_W(V\oplus \R)$ in each fiber.

\begin{prop}\label{prop:codim zero cobordism}
Let $(M,f)$ be a $V$-framed cobordism element of $(X,A)$ and $N\subset M$ a $G$-invariant codimension $0$ submanifold of the interior of $M$ such that $f(\overline{M\setminus N})\subset A$. Then $(M,f)$ and $(N,f|_N)$ are $V$-framed cobordant. 
\end{prop}

\begin{proof}
Define the space $W$ to be the pushout
\begin{figure}[H]
\center
\begin{tikzcd}
N\times\{0\} \arrow[r]\arrow[d] & N\times[0,1]\arrow[d]\\
M\times[-1,0]\arrow[r]&W \arrow[ul, phantom, "\ulcorner", very near start]
\end{tikzcd}
\end{figure}
While we would like to say that $W$ is a cobordism from $M$ to $N$, this is not technically true. Firstly, $W$ is not a manifold with corners since near $\partial N\times\{0\}$, $W$ is locally homeomorphic to $\R^{n-1}\times[0,\infty)\cup\R^{n-1}\times[0,\infty)$. Secondly, $W$ has corners at $\partial M\times\{0\}$ that do not coincide with $\partial M\times\{-1\}\cup \partial N\times\{1\}$. We will replace $W$ with another space, $W'$, such that $W'$ is homeomorphic to $W$ and is a cobordism from $M$ to $N$.

We will abuse notation and refer to the part of $W$ not locally homeomorphic to $\R^n$ as $\partial W$ and the part of $W$ that is locally homeomorphic to $\R^n$ as the interior of $W$. Begin by choosing an inward pointing vector field $\xi$ on $W$. By an inward pointing vector field, we mean that for each $x\in \partial W$, the vector at $x$ points to the interior of whichever kind of chart is providing the local model at $x$. Such a vector field can be built by gluing together local vector fields with a smooth partition of unity. By averaging across $G$, we may assume that $\xi$ is $G$-invariant. We put a smooth structure on $\partial W$ as in \cite[Section 6]{manifold_approach} by taking disks in the interior of $W$ transverse to $\xi$ and flowing backward along $\xi$ to reach $\partial W$.

Since $M$ is compact and $\xi$ is inward pointing, we have a well defined flow map associated to $\xi$ for all positive times. This gives rise to a map $\phi:\partial W\times [0,\infty)\to W$. However, $\phi$ is not smooth since we have changed the smooth structure on $\partial W$. Define a bump function
\begin{align*}
\Psi:\R&\to \R\\
t & \mapsto \begin{cases}
\exp\left(\tfrac{1}{(2t)^2-1}\right) & \text{ if $t\in (-1/2,1/2)$}\\
0 & \text{ else }
\end{cases}.
\end{align*} Now define $U=\phi(\partial W\times (0,\Psi(t)))$ with smooth structure arising from the fact that it is an open subset of the interior of $W$. The projection map $p:U\to \partial W$ is smooth by construction. In fact, wherever $U$ is non-empty, $p$ is a smooth submersion with fibers homeomorphic to open intervals. Therefore, wherever $U$ is non-empty, $p$ has a smooth section, $s$. Note also that by construction, $s$ continuously approaches $\partial W$ as $t$ approaches both $-1/2$ and $1/2$. In order to make $s$ smoothly approach $\partial W$ we multiply by a bump function near both $1/2$ and $-1/2$. Then the image of this smoothing of $s$ along with 
\[
M\times\{-1\}\cup\partial M\times[-1,-1/2] \amalg N\times\{1\}\cup \partial N\times[1/2,1]
\]
defines the boundary of $W'$ a $G$-cobordism from $M$ to $N$. Moreover, the inclusion $W'\to W$ is homotopic to a homeomorphism. 

So far we have built a $G$-cobordism, $W'$, from $M$ to $N$. To see that $M$ and $N$ are $V$-framed cobordant we first observe that the $V$-framing of $N$ arises from restricting that of $M$. Thus, $M\times I$ and $N\times I$ are $(V\oplus\R)$-framed in a compatible way. Then $W'$ inherits a $(V\oplus\R)$-framing from those on $M\times I$ and $M\times I$. This framing on $W'$ immediately recovers the $V$-framing on $M\amalg N$. The map $W'\to X$ recovering $f$ and $f|_N$ is inherited from the corresponding maps on $M\times I$ and $N\times I$ as follows. The inclusion $W'\hookrightarrow W$ is continuously homotopic, relative to $M\times\{-1\}\amalg N\times\{1\}$, to a homeomorphism. So we may build a map $W'\to X$ by composing $W\to X$ with this homeomorphism.  This new map agrees with $f$ on $M\times\{-1\}$ and $N\times\{1\}$ while mapping the sides of $W'$ into $A$ as desired.
\end{proof}

\begin{df}\label{df:neat}
A closed embedding of smooth manifolds with corners is \textbf{neat} if it is locally modelled on the inclusion
\[
[0,\infty)^2\times \R^{m-2}\times \{0\}^{n-m} \hookrightarrow [0,\infty)^2\times \R^{m-2}\times \R^{n-m}
\]
\end{df}

\begin{rmk}
A neat embedding $i:M\to N$ of manifolds with boundary will always have the property that $i(\partial M)\subset \partial N$. 
\end{rmk}

\begin{rmk}
If $M$ is a smooth $G$-manifold with corners, as in \cref{df:G-manifold with corners}, then $M^H$ is a neat submanifold of $M$ for each $H\leq G$ where neatness is defined in \cref{df:neat}.
\end{rmk}

On occasion, we will work with the normal bundle of a cobordism, $W$, from $M$ to $N$. In order for $W$ to have a well-defined normal bundle, we must carefully choose an appropriate embedding of $W$ into a larger manifold. Let $U$ be some finite dimensional $G$-representation such that $W$ embeds neatly into $U\times [0,1]\times \R_{\geq 0}$. We then choose our embedding so that:
\begin{enumerate}
\item $M$ embeds in $U\times \{0\}\times \R_{\geq 0}$ and
\item $N$ embeds in $U\times \{1\}\times \R_{\geq 0}$  
\end{enumerate}
This embedding guarantees that the sides of $W$ embed into $U\times[0,1]\times\{0\}$. Moreover, the codimension of each point of $W$ in $U\times[0,1]\times\R_{\geq 0}$ agrees and we can construct a well defined normal bundle. 

\begin{rmk}
As a consequence of requiring $W\hookrightarrow U\times [0,1]\times \R_{\geq 0}$ to be neat, we are guaranteed that $\partial^cW$ embeds in $U\times\{0,1\}\times \{0\}$.
\end{rmk}

\begin{lem}\label{lem:eq cobordism is eq rel}
Equivariant framed cobordism is an equivalence relation.
\end{lem}

\begin{proof}
Let $(M,f)$ be a $V$-framed cobordism element of $(X,A)$. Then $M\times I$ is a cobordism from $M$ to $M$ (where $G$ acts trivially on $I$). Since the tangent bundle commutes with products, $M\times I$ is $(V\oplus \R)$-framed and induces the given framing on the top and bottom. The reference map $q:M\times I\to X$ is given by $q(m,t)=f(m)$ which certainly restricts to $f$ on the boundary of $M\times I$. 

Any cobordism from $M$ to $N$ is also a cobordism from $N$ to $M$. 

Now let $(W,q)$ be a cobordism from $(M,f)$ to $(M',f')$ and $(W',q')$ a cobordism from $(M',f')$ to $(M'',f'')$. We choose collared neighborhoods $M'\times (-1,0]$ in $W$ and $M'\times [0,1)$ in $W'$. The existence of these collared neighborhoods is shown in \cite{gimm}. Then $M'\times(-1,0]$ and $M'\times[0,1)$ are equivariantly diffeomorphic as framed manifolds. We can then glue together the collared neighborhoods along $M'\times(-\epsilon,\epsilon)$ to produce a smooth framed $G$-manifold with bottom face $M$ and top face $M''$. To build the reference map to $X$ on this new manifold we observe that $q|_{M'\times I}$ and $q'|_{M'\times I}$ are homotopic. We then create a new map by applying $q$ on $W$, homotoping $q$ to $q'$ where we glue $W$ and $W'$ on the overlap, and then applying $q'$ on $W'$.
\end{proof}

\begin{lem}
The set of $V$-framed cobordism elments of $(X,A)$ is an abelian group under disjoint union of $G$-manifolds, which we denote $\omega_V^G(X,A)$.
\end{lem}
\begin{proof}
Since $\omega_V^G(X,A)$ is a commutative monoid with identity element given by $\emptyset$, it suffices to show that the elements of $\omega_V^G(X,A)$ have inverses with respect to disjoint union. 

Let $(M,f)\in \omega_V^G(X,A)$. We claim that $(-M,f)$ is the inverse of $(M,f)$. Let 
\[
\phi:TM\oplus\varepsilon_M(\R^n)\cong \varepsilon_M(V)\oplus\varepsilon_M(\R^n)
\]
be the framing of $M$ and $-\phi$ be the corresponding framing for $(-M)$. Let $N=M\times I$ where $G$ acts trivially on $I$. Then $N$ is a $G$-manifold. To see that $N$ is $(V\oplus \R)$-framed, observe that

\begin{align*}
TN\oplus\varepsilon_N(\R^{n+1})
& \cong TM\oplus TI \oplus\varepsilon_M(\R^n)\oplus\varepsilon_I(\R)\\
& \cong TM\oplus\varepsilon_M(\R^n)\oplus\varepsilon_I(\R)\oplus\varepsilon_I(\R)\\
& \cong \varepsilon_M(V)\oplus \varepsilon_M(\R^n)\oplus\varepsilon_I(\R)\oplus\varepsilon_I(\R)\\
& \cong \varepsilon_{M\times I}(V\oplus \R) \oplus \varepsilon_{M\times I}(\R^n\oplus \R)\\
& \cong \varepsilon_N(V\oplus \R) \oplus \varepsilon_N(\R^{n+1})
\end{align*}

Moreover, $M\amalg (-M)\subset \partial N$, and the depth two corners of $N$ coincide with $\partial M \amalg \partial (-M)$. The $V$-framing induced by restriction of the framing on $N$ to $M\amalg (-M)$ agrees with the $V$-framing on $M\amalg (-M)$. Lastly, we see that $q=(f,1_I):M\times I=N\to X$ is an equivariant map with $q|_M=f$ and $q|_{(-M)} =f$ and $q$ maps the sides of $N$ into $A$. 
\end{proof}

\begin{rmk}
The groups $\omega_V^G(X,A)$ also carry the structure of a Mackey functor with $\res^G_H$ given by forgetting a $G$-action to an $H$-action and $\tr_H^G$ given by $G\times_H(-)$. However, we will not need this structure in the present paper.
\end{rmk}

\begin{rmk}
For each pair of finite dimensional orthogonal $G$-representations, $V$ and $W$, there is a suspension isomorphism $\omega_V^G(X,*)\to \omega^G_{V\oplus W}(\Sigma^W X,*)$ given by taking the product with $D(W)$ and then trimming the corners as in \cite{gimm}. The functors which take a $G$-space, $X$, and a $G$-representation, $V$, to $\omega_V^G(X,*)$ along with the suspension isomorphisms form an $RO(G)$-graded homology theory in the precise sense of \cite[Chapter XIII]{alaska}.
\end{rmk}

\beforesubsection
\subsubsection{Equivariant stable homotopy theory}\aftersubsection

We will now briefly introduce some preliminaries on equivariant stable homotopy theory. Further background on equivariant stable homotopy theory can be found in \cite{mandell2002equivariant}, \cite{schwede_global}, and \cite{hhr_book}.

\begin{df}
We denote the \textbf{category of orthogonal $G$-spectra} as $G\Osp$. While we will primarily work with orthogonal $G$-spectra, it will, on occasion, be more convenient to work with diagram spectra. There is an equivalence of categories between $G\Osp$ and the category of functors $\textbf{J}_G\to \GTop_*$ enriched in $\GTop_*$ \cite{mandell2002equivariant}. Here $\textbf{J}_G$ is a category of finite dimensional $G$-representations. We will refer to this category of diagram spectra as \textbf{genuine $G$-spectra}.
\end{df}

\begin{df}
Let $\rho$ be the regular $G$-representation, let $V$ be a finite dimensional orthogonal $G$-representation, and let $X$ be a genuine $G$-spectrum. Define the \textbf{$V$th equivariant stable homotopy group} of $X$ as 
\[
\pi_V^G(X) = \underset{n\to\infty}{\colim}[S^{n\rho\oplus V},X(n\rho)]_*^G
\]
where the colimit on the right is of based $G$-homotopy classes of continuous $G$-maps.
\end{df}

\begin{rmk}
The category of orthogonal $G$-spectra has a model structure developed in \cite[III.4.2]{mandell2002equivariant}. Throughout this section, when we derive functors on orthogonal $G$-spectra we are using the (co)fibrant replacement from this model structure unless otherwise noted. 
\end{rmk}

\begin{df}\label{df:genuine and categorical fixed points}
Let $X\in G\Osp$. The \textbf{categorical fixed points} of $X$ is the orthogonal $WH$-spectrum with $n$th level given by $X_n^H$. The \textbf{genuine fixed points} of $X$, denoted $X^G$, are the right derived categorical fixed points. This means that $X^G$ is constructed by replacing $X$ with a stably equivalent fibrant spectrum and then taking the categorical fixed points.
\end{df}

Note that the $G$-equivariant stable homotopy groups of $X$ agree with the non-equivariant stable homotopy groups of $X^G$. For details on various fixed point functors see \cite{hhr_book} or \cite{mandell2002equivariant}.

\beforesubsection
\subsection{Isotropy separation}\aftersubsection

In this section we will define four maps and then show that they fit into a commuting square. First we need a lemma relating trivializations of tangent and normal bundles.

\begin{lem}\label{lem:normal bundle}
Let $G$ be a compact Lie group and $(M,f)$ a $V$-framed cobordism element of $(X,A)$. There exists a finite dimensional real orthogonal $G$-representation $W$ such that the normal bundle of $M$ embedded in $kW\oplus V$ is isomorphic to $M\times kW$ for sufficiently large $k$. 
\end{lem}

\begin{proof}
Let $(M,f)$ be a $V$-framed cobordism element of $(X,A)$ so that
\begin{equation}\label{eq:eq1}
TM\oplus \varepsilon_M(\R^n)\cong \varepsilon_M(V)\oplus\varepsilon_M(\R^n).
\end{equation}
There exists some finite dimensional orthogonal $G$-representation $W$ and large $k$ so that we can embed $M$ into $(k-n)W$ as in \cite[Lemma 2.21]{gimm}. Letting $\nu(M,(k-n)W)$ be the normal bundle of this embedding, we have the following bundle isomorphism:

\begin{equation}\label{eq:eq2}
TM\oplus\nu(M,(k-n)W) \cong \varepsilon_M((k-n)W).
\end{equation}
Combining \cref{eq:eq1} and \cref{eq:eq2} yields:
\[
\varepsilon_M(V)\oplus\nu(M,(k-n)W)\oplus\varepsilon_M(\R^n) \cong \varepsilon_M((k-n)W)\oplus\varepsilon_M(\R^n).
\]
Therefore, 
\[
\nu(M,(k-n)W\oplus \R^n\oplus V)\cong\varepsilon_M(\R^n)\oplus\varepsilon_M((k-n)W).
\]

We can thus find $k$ large enough so that
\begin{equation}\label{eq:trivialize normal bundle}
\nu(M,kW\oplus V) \cong \varepsilon_M(kW).
\end{equation}
\end{proof}

\begin{df}\label{df:equivariant framed PT}
Let $G$ be a compact Lie group. The \textbf{equivariant framed Pontryagin-Thom construction}, denoted $\Delta$, is a map
\[
\omega_V^G(X,A) \xrightarrow{\Delta} \pi_V^G\left( \Sigma^\infty X/A \right)
\]
defined as follows. Let $(M,f)$ be a $V$-framed cobordism class of $(X,A)$. By \cref{lem:normal bundle} we may choose a finite dimensional orthogonal $G$-representation, $W$, and an integer, $k$, so that $\nu(M,kW\oplus V)\cong M\times kW$. Using this large $k$, embed $M$ into $S^{kW\oplus V}$ so that it is far from the basepoint. Define $\Delta(M,f)$ to be the following composite of maps:
\begin{align*}
S^{kW\oplus V} 
& \xrightarrow{\cong} D(kW\oplus V)/S(kW\oplus V)\\
& \xrightarrow{} D(\nu)/(D(\nu|_{\partial M})\cup S(\nu))\\
& \xrightarrow{\cong} (M\times D(kW))/((\partial M\times D(kW))\cup (M\times S(kW)))\\
& \xrightarrow{f\times 1} (X\times D(kW))/((A\times D(kW)) \cup (X\times S(kW))).
\end{align*}
The second map in the composite is a Pontryagin-Thom collapse map. In other words, the second map sends a point in the interior of the disc bundle of $\nu$ to itself and a point outside the interior of $D(\nu)$ to the basepoint. Note that in order to guarantee that $D(\nu)\subset D(kW\oplus V)$ we may have to rescale the unit disc to be sufficiently large. Now identifying the target of the above composition with $S^{kW}\wedge X/A$ we obtain an element of $\pi_V^G \left( \Sigma^\infty X/A\right)$. 
\end{df}

\begin{lem}
The map $\Delta$ of \cref{df:equivariant framed PT} is well defined. 
\end{lem}
\begin{proof}
Let $(M,f)$ and $(M',f')$ be two $V$-framed $G$-manifolds which are cobordant along $(N,q)$. Embed $N$ into $S^{kW\oplus V}\times I$ so that $M$ and $M'$ embed into $S^{kW\oplus V}\times\{0\}$ and $S^{kW\oplus V}\times\{1\}$ respectively. Then applying the rest of the construction of \ref{df:equivariant framed PT} to $N\hookrightarrow S^{kW\oplus V}\times I$ gives an equivariant homotopy from $\Delta(M,f)$ to $\Delta(M',f')$.
\end{proof} 

\begin{rmk}
We will now restrict our attention to the case where $G$ is isomorphic to the product of a finite group and a torus. As noted in \cite[6.2.13]{schwede_global}, this is the largest class of groups for which the map of \cref{df:equivariant framed PT} can be expected to be an isomorphism. Schwede shows that the composition of induction and the Wirthm\"uller map in equivariant cobordism
\[
\Omega_V^H(X) \to \Omega_{V\oplus L}^H\left(X_+\wedge S^L,*\right)
\] 
is given by exterior multiplication by the cobordism class of $S(\R\oplus L)$ where $L$  is the $H$-representation $T_{eH}G/H$ \cite[6.2.12]{schwede_global}. Since induction is an isomorphism, the Wirthm\"uller map is invertible precisely when exterior multiplication by $S(\R\oplus L)$ is invertible. Here invertibility fails if $H$ acts non-trivially on $L$. When $G$ is isomorphic to a product of a finite group and a torus, this action is trivial. 
\end{rmk}

\begin{df}
Given a $G$-space $X$, we define the $H$-fixed points of $X$ to be
\[
X^H=\{x\in X\mid hx=x \text{ for all } h\in H\}.
\]
Given an equivariant map between $G$-spaces, $f:X\to Y$ we define $f^H:X^H\to Y^H$ to be the restriction of $f$ to $X^H$. Note that this restriction lands in $Y^H$ since $hf(x)=f(hx)=f(x)$ for all $x\in X^H$ and $h\in H$. 
\end{df}

\begin{df}\label{df:fixed points of cobordism}
The \textbf{$H$-fixed points of a cobordism class} is a map
\[
\omega_V^G(X,A) \xrightarrow{\Psi_\omega^H} \omega_{V^H}^{WH}(X^H,A^H)
\]
given by sending $(M,f)$ to $(M^H,f^H)$. Note that $M^H$ is a $V^H$-framed $WH$-manifold and that $f^H:(M^H,\partial M^H)\to (X^H,A^H)$ is $WH$-equivariant. Here $WH=N_GH/H$ is the Weyl group of $H$.
\end{df}

\begin{df}\label{df:fixed points homotopy theory}
We define
\[
\pi_V^G \left( \Sigma^\infty X/A \right) \xrightarrow{\Psi_\pi^H} \pi_{V^H}^{WH}\left( \Sigma^\infty X^H/A^H \right)
\]
to be the restriction map from the genuine $G$-fixed points to genuine $WH$-fixed points of the geometric $H$-fixed points. In other words, we send the $G$-homotopy class of the $G$-equivariant map 
\[
S^{k\rho\oplus V}\rightarrow S^{k\rho}\wedge X/A
\]
to the $WH$-homotopy class of the $WH$-equivariant map 
\[
S^{k\rho^H\oplus V^H}\rightarrow S^{k\rho^H}\wedge X^H/A^H
\]
\end{df}

\begin{df}\label{df:fixed points PT}
The \textbf{equivariant framed Pontryagin-Thom construction on fixed points} is a map 
\[
\omega_{V^H}^{WH}(X^H,A^H)\xrightarrow{\Delta^H} \pi_{V^H}^{WH}\left( \Sigma^\infty X^H/A^H \right)
\] 
given by taking the $H$-fixed points of everything in \cref{df:equivariant framed PT}.
\end{df}

\begin{prop}\label{prop:big diagram commutes}
The following diagram commutes:

\begin{diagram}\label{dia:big diagram}
\omega_V^G(X,A) \arrow[r,"\Delta"] \arrow[d,"\Psi_\omega^H"] & \pi_V^{G}(\Sigma^\infty X/ A)\arrow[d,"\Psi_\pi^H"]\\
\omega_{V^H}^{WH}(X^H,A^H) \arrow[r,"\Delta^H"] & \pi_{V^H}^{WH}(\Sigma^\infty X^H/A^H)
\end{diagram}

where the maps $\Delta$, $\Psi_\omega^H$, $\Psi_\pi^H$ and $\Delta^H$ are the maps of \cref{df:equivariant framed PT}, \cref{df:fixed points of cobordism}, \cref{df:fixed points homotopy theory}, and \cref{df:fixed points PT} respectively. 
\end{prop}

\begin{proof}
The proof follows from the fact that taking the $H$-fixed points of the Pontryagin-Thom map, is the same as applying the Pontryagin-Thom map to the $H$-fixed points. 
\end{proof}

\beforesubsection
\subsection{Proof of the equivariant Pontryagin-Thom theorem}\aftersubsection

In this section we will prove our first main theorem.

\begin{thm}
Let $G$ be isomorphic to a product of a finite group and a torus. There is an isomorphism of abelian groups:
\[
\omega_V^G(X,A)\cong \pi_V^G(\Sigma^\infty X/A).
\]
\end{thm}

We will prove this using an isotropy separation technique. More specifically, we will first prove the result for free group actions and then piece various isomorphisms together using equivariant classifying spaces. 

\begin{df}
Let $G$ be a group. A \textbf{family of subgroups of $G$} is a set of subgroups of $G$ closed under both conjugation and taking subgroups.
\end{df}

\begin{df}
Let $G$ be a group, and $F$ be a family of subgroups of $G$. The \textbf{universal space associated to the family $F$}, denoted $EF$, is the terminal object in the category of $G$-spaces with isotropy in $F$. It is defined, up to $G$-weak equivalence, by the formula
\[
EF^H = \begin{cases}
* & \text{ if $H\in F$}\\
\emptyset & \text{ if $H\notin F$}.
\end{cases}
\]
\end{df}

\begin{notn}
Let $F$ and $F'$ be families of subgroups such that $F'\subset F$. We define
\[
\omega_V^G[F,F'](X,A) = \omega_V^G(EF\times X, EF\times A\cup EF'\times X)
\]
and
\[
\pi_V^G[F,F'](X,A) = \pi_V^G(\Sigma^\infty EF\times X/EF\times A\cup EF'\times X).
\]

Note that $EF\times X/EF\times A\cup EF'\times X\cong EF/EF'\wedge X/A$.

If $F'=\emptyset$ and $F=\{1\}$ or $F=\{\All\}$ we write
\[
h_V^G[\free](X,A) \qquad \text{or} \qquad h_V^G[\All](X,A)
\]
respectively for $h=\omega,\pi$. This notation is aligned with the fact that $E\emptyset\simeq \emptyset$, and $E\{1\}\simeq EG$, and $E\{\All\}\simeq *$. 
\end{notn} 

\begin{notn}
Let $G$ be a group, and $H\leq G$. Denote the conjugacy class of $H$ as $(H)$. 
\end{notn}

Now let $F$ and $F'$ be two adjacent families of subgroups so that $F=F'\cup (H)$. By Elmendorf's theorem, the $H$-fixed points of $EF\times X$ is $EWH\times X^H$. Thus, replacing $(X,A)$ with $(EF\times X,EF\times A\cup EF'\times X)$ in \cref{dia:big diagram} we immediately see that the following diagram commutes:

\begin{diagram}\label{dia:isotropy separation diagram}
\omega_V^G[F,F'](X,A)\arrow[d,"\Psi_\omega^H"]\arrow[r,"\Delta"]& \pi_V^{G}[F,F'](\Sigma^\infty X/ A)\arrow[d,"\Psi_\pi^H"]\\
\omega_{V^H}^{WH}[\free](X^H,A^H)\arrow[r,"\Delta^H"]&\pi_{V^H}^{WH}[\free](\Sigma^\infty X^H/ A^H).
\end{diagram}

\begin{prop}\label{prop:free PT isomorphism}
$\Delta^H$ is an isomorphism. 
\end{prop}
\begin{proof}
Cobordism elements in $\omega_{V^H}^{WH}[\free](X,A)$ are free $WH$-manifolds. The proof that $\Delta^H$ is an isomorphism then follows from the usual proof of the Pontryagin-Thom theorem (see \cite{stong68}) leveraging the fact that an equivariant map between manifolds with free actions may be replaced by a $WH$-homotopic map which is transverse to a submanifold of our choosing. 
\end{proof}

Our next goal is to show that the two vertical maps in the above diagram are isomorphisms.

\subsubsection{The map $\Psi_\omega^H$ is an isomorphism}

If $H$ is a subgroup of $G$, then we denote the normalizer of $H$ in $G$ by $NH$. 

Recall that $G\times_{NH}(-)$ is left adjoint to the functor that forgets a $G$-action down to an $NH$-action. The counit, $G\times_{NH}X\to X$, is given by $(g,x)\mapsto gx$. The unit, $X\to G\times_{NH}X$ is given by $x\mapsto (e,x)$.  

\begin{df}\label{df:the manifold Q}
Define a map
\[
\omega_{V^H}^{WH}[\free](X^H,A^H) \xrightarrow{\Theta} \omega_V^G[F,F'](X,A)
\]
by $\Theta(N,t)=(Q,q)$ where $Q$ is given by
\[
Q = G\times_{NH}(N\times D(V^{H^\perp})).
\]

Here $V^{H^\perp}$ is the orthogonal complement of $V^H$ in $V$. The map $q:Q\to EF\times X$ is constructed as follows. Since all the isotropy groups of $Q$ are contained in $F$, there is a map 

\[
q_1:Q\to EF.
\]

Let $p:N\times D(V^{H^\perp})\to N$ be the projection. Then define the map $q:Q\to EF\times X$ as

\[
G\times_{NH}(N\times D(V^{H^\perp})) \xrightarrow{q_1\times G\times_{NH}(p)} EF\times G\times_{NH} N \xrightarrow{\id\times G\times_{NH}(t)} EF\times G\times_{NH} X^H \hookrightarrow EF\times X.
\]
\end{df}

\begin{rmk}
Recall that $G$ is a product of a finite group and a torus, and write $G\cong F\times T$ where $F$ is finite and $T$ a torus. The normalizer of $H\leq G$ must have the same dimension as $G$. This follows immediately from the fact that $\{e\}\times T\leq NH$ for all $H\leq G$. Note, that the dimension of $Q$ in the previous definition is $\dim (V)$. 
\end{rmk}

Unfortunately, $Q$ is a manifold with corners whereas we wish for it to be a manifold with boundary. We rectify this by replacing $Q$ with a trimming of $Q$ as in \cite[2.22]{gimm}. This trimming is given by cutting off a collar of the boundary of $Q$ in such a way that the remaining manifold is a manifold with boundary rather than with corners. This new manifold (which we also call $Q$) is equivariantly homeomorphic to $Q$. We obtain the map from the trimmed version of $Q$ to $EF\times X$ simply by restricting the original map on $Q$ to the trimming. 

\begin{lem}
The manifold 
\[
Q = G\times_{NH}(N\times D(V^{H^\perp}))
\]
is a $V$-framed $G$-manifold.
\end{lem}
\begin{proof}
The $G$-action on $Q$ is given by first pulling back the $WH$-action on $N$ along the quotient $NH\rightarrow WH$ and then taking the induced action.

Let
\[
\phi:TN\oplus\varepsilon_N(\R^k)\xrightarrow{\cong} \varepsilon_N(V^H\oplus\R^k)
\]
be the $V^H$-framing of $N$.

Observe that
\begin{align*}
T(N\times D(V^{H^\perp})) \oplus \varepsilon_{N\times D(V^{H^\perp})}(\R^k)
& \cong [TN \oplus \varepsilon_N(\R^k)]\times[TD(V^{H^\perp})]\\
& \cong \varepsilon_N(V^H\oplus \R^k)\times (D(V^{H^\perp})\times V^{H^\perp})\\
& \cong \varepsilon_{N\times D(V^{H^\perp})}(V\times \R^k)
\end{align*}
Therefore, $N\times D(V^{H^\perp})$ inherits a $V$-framing from the $V^H$-framing on $N$. Then $Q$ is $V$-framed as well since 
\[
TQ = T(G\times_{NH}(N\times D(V^{H^\perp})))\cong G\times_{NH} T(N\times D(V^{H^\perp}))
\]
Here we are using the fact that $G\times_{NH}(-)$ commutes with the tangent bundle since $NH$ is codimension zero in $G$. 

Note that when we replace $Q$ with a trimming, the trimming is also $V$-framed  by extending the $V$-framing on $Q$ along the homeomorphism from $Q$ to its trimming. 
\end{proof}

\begin{lem}
The image of the restriction of $q:Q\to EF\times X$ to $\partial Q$ lies in $(EF\times A)\cup (EF'\times X)$.
\end{lem}

\begin{proof}
First observe that 
\[
\partial Q = G\times_{NH}(\partial N\times D(V^{H^\perp}))\cup G\times_{NH}(N\times S(V^{H^\perp}))
\]
Also note that the isotropy groups of $N\times D(V^{H^\perp})$ are always contained in $H$ and that this containment is proper everywhere except $N\times\{0\}$. In particular, the isotropy groups of $N\times S(V^{H^\perp})$ are all contained in $F'$. We can therefore see that $q:Q\to EF\times X$ restricts to
\[
G\times_{NH}(N\times S(V^{H^\perp}))\rightarrow EF'\times X
\]
and
\[
G\times_{NH}(\partial N \times D(V^{H^\perp}))\to EF\times A
\]
Therefore,
\[
q|_{\partial Q}:G\times_{NH}(\partial N\times D(V^{H^\perp}))\cup G\times_{NH}(N\times S(V^{H^\perp})) \rightarrow EF'\times X\cup EF\times A
\]
as desired. 
\end{proof}

The previous two lemmas show that $(Q,q)$ is indeed an element of $\omega_V^G[F,F'](X,A)$. 

\begin{lem}
The map 
\[
\Theta:\omega_{V^H}^{WH}[\free](X^H,A^H) \rightarrow \omega_V^G[F,F'](X,A)
\]
is well defined.
\end{lem}

\begin{proof}
Let $(N_1,t_1)=(N_2,t_2)$ in $\omega_{V^H}^{WH}[\free](X^H,A^H)$ and say $(W,w)$ is a $(V^H\oplus \R)$-framed $WH$-cobordism between them. Then passing $(W,w)$ through $\Theta$ will give a $(V\oplus \R)$-framed $G$-cobordism between $\Theta(N_1,t_1)$ and $\Theta(N_2,t_2)$. 
\end{proof}

\begin{prop}\label{prop:right inverse}
The map $\Theta$ is a right inverse to $\Psi_\omega$. In other words
\[
\Psi_\omega \circ \Theta = \id
\]
\end{prop}

\begin{proof}
Let $(N,t)\in \omega_{V^H}^{WH}[\free](X^H,A^H)$ and let $(Q,q)=\Theta(N,t)$. Now observe that 
\[
Q^H = (G\times_{NH}(N\times D(V^{H^\perp})))=N^H=N
\]
where the last equality follows from the fact that $N$ is a $WH$-manifold. This also shows that we recover the original framing on $N$.

Also recall from \cref{df:the manifold Q} that 
\[
q=q_1\times (G\times_{NH}(t\circ p)): Q\rightarrow EF\times (G\times_{NH} X^H)\rightarrow EF\times X
\]
Thus,
\[
q^H = q_1^H \times (G\times_{NH}(t\circ p)^H) = t^H=t
\]
Here $q_1^H$ is a constant map since $H\in F$ and $t^H=t$ since $N$ is a free $WH$-manifold. 

Therefore, $\Psi_\omega^H\circ \Theta = \id$ as desired. 
\end{proof}

Before showing that $\Theta$ is a left inverse to $\Psi_\omega^H$, we must first prove two technical lemmas. 

\begin{lem}\label{lem:tubular nbhd}
Let $M$ be a $V$-framed $G$-manifold. Let $\nu(M^H,M)$ denote the normal bundle of $M^H\hookrightarrow M$. Then $\nu(M^H,M)\cong M^H\times V^{H^\perp}$.
\end{lem}

\begin{proof}
Let $M$ be a $V$-framed $G$-manifold with framing
\begin{equation}\label{eq:trivialization}
TM\oplus\varepsilon_M(\R^n)\cong \varepsilon_M(V)\oplus\varepsilon_M(\R^n)
\end{equation}
Now observe that 
\begin{align*}
TM^H\oplus \nu(M^H,M)\oplus \varepsilon_{M^H}(\R^n)
& \cong (TM\oplus\varepsilon_M(\R^n))|_{M^H}\\
& \cong (\varepsilon_M(V)\oplus\varepsilon_M(\R^n))|_{M^H}\\
& \cong \varepsilon_{M^H}(V^H)\oplus \varepsilon_{M^H}(V^{H^\perp})\oplus\varepsilon_{M^H}(\R^n)
\end{align*}
where the second equality follows from \cref{eq:trivialization}. In short, we obtain
\begin{equation}\label{eq:normal bundle}
TM^H\oplus \nu(M^H,M)\oplus \varepsilon_{M^H}(\R^n) \cong \varepsilon_{M^H}(V^H)\oplus \varepsilon_{M^H}(V^{H^\perp})\oplus\varepsilon_{M^H}(\R^n)
\end{equation}

Each of the summands in \cref{eq:normal bundle} is a bundle over $M^H$. Therefore, the action of $H$ on each of these bundles is fiberwise. We may thus take the orthogonal complement of the $H$-fixed points in each fiber of the bundles in \cref{eq:normal bundle} while preserving the isomorphism. 

First, note that $H$ acts trivially on $TM^H$ and $\varepsilon_{M^H}(\R^n)$ and $\varepsilon_{M^H}(V^H)$. Therefore, the fiberwise orthogonal complements of the $H$-fixed points in each of these three bundles will be a zero dimensional bundle over $M^H$. 

Second, observe that the action of $H$ on the non-zero vectors in $\varepsilon_{M^H}(V^{H^\perp})$ and $\nu(M^H,M)$ is fixed point free. Therefore, taking the fiberwise orthogonal complement of the $H$-fixed points in these two bundles, does not change the bundle at all. 

In conclusion, taking fiberwise orthogonal complements of the $H$-fixed points of the bundles in \cref{eq:normal bundle} yields
\[
\nu(M^H,M)\cong \varepsilon_{M^H}(V^{H^\perp})
\]
as desired. 
\end{proof}

\begin{lem}\label{lem:embedding}
There is a $G$-equivariant embedding $G\times_{NH}(M^H\times D(V^{H^\perp}))\hookrightarrow M$ which lands in the interior of $M$.
\end{lem}

\begin{proof}
By \cref{lem:tubular nbhd}, and the tubular neighborhood theorem for $G$-manifolds with boundary, $M^H\hookrightarrow M$ has a tubular neighborhood diffeomorphic to $M^H\times D(V^{H^\perp})$. We  claim that the induced map 
\[
G\times_{NH}M^H \to M
\]
is injective. Suppose that this map fails to be injective. Then there is some $x\in M$ such that $x\in M^H\cap M^{gHg^{-1}}$ with $H\neq gHg^{-1}$. Then both $H$ and $gHg^{-1}$ are contained in $G_x$. Moreover, since $H\neq gHg^{-1}$ the containment of $H$ and $gHg^{-1}$ in $G_x$ must be proper. Recall however that the isotropy of each point in $G\times_{NH}(M^H\times D(V^{H^\perp}))$ is in $F$. By construction if $G_x\in F$ and $H\leq G_x$ then $H=G_x$ giving us a contradiction.

Since we have injections $M^H\times D(V^{H^\perp})$ and $G\times_{NH} M^H$, we may obtain an injection $G\times_{NH}(M^H\times D(V^{H^\perp}))$ by shrinking the radius of $D(V^{H^\perp})$ as necessary. The equivariant tubular neighborhood theorem guarantees that this injection is $G$-equivariant.

So far we have shown that $G\times_{NH}(M^H\times D(V^{H^\perp}))$ embeds equivariantly into $M$. To see that we can embed $G\times_{NH}(M^H\times D(V^{H^\perp}))$ into the interior of $M$ observe that if the inclusion coincides with $\partial M$ we may choose a $G$-invariant inward pointing vector field and flow along it to perturb the image of $G\times_{NH}(M^H\times D(V^{H^\perp}))$ off of the boundary while preserving equivariance. 
\end{proof}

\begin{prop}\label{prop:left inverse}
The map $\Theta$ is a left inverse for $\Psi_\omega^H$. In other words:
\[
\Theta\circ \Psi_\omega^H = \id
\]
\end{prop}

\begin{proof}
First, let 
\[
j:G\times_{NH}(M^H\times D(V^{H^\perp})) \to M
\]
be the $G$-equivariant embedding into the interior of $M$ from \cref{lem:embedding}. 

By construction, the map
\[
q:G\times_{NH}(M^H\times D(V^{H^\perp})) \rightarrow X
\]
is the restriction of $f:M\rightarrow X$ to $G\times_{NH}(M^H\times D(V^{H^\perp}))$.

Recall that the normal bundle of $M^H$ in $M$ is $M^H\times V^{H^\perp}$. Therefore, setting 
\[
Q=G\times_{NH}(M^H\times D(V^{H^\perp})),
\]
the following diagram commutes and each of the maps is an isomorphism:

\begin{diagram*}
TQ\oplus\varepsilon_Q(\R^n) \arrow[d,"\cong"] \arrow[r,"\cong"] & \varepsilon_Q(V^H)\oplus \varepsilon_Q(V^{H^\perp})\oplus\varepsilon_Q(\R^n) \arrow[d,"\cong"]\\
(TM\oplus\varepsilon_M(\R^n))|_Q \arrow[r,"\cong"]&(\varepsilon_M(V)\oplus\varepsilon_M(\R^n))|_Q.
\end{diagram*}

Thus, the tangent bundle of $G\times_{NH}(M^H\times D(V^{H^\perp}))$ agrees with the tangent bundle of $M$ restricted to $G\times_{NH} (M^H\times D(V^{H^\perp}))$. 

To summarize, the map $j:Q=G\times_{NH}(M^H\times D(V^{H^\perp}))\to M$ embeds $Q$ into the interior of $M$ in such a way so that the restriction of $f:M\to X$ to $Q$ gives $q$ and the restriction of the framing of $M$ to $Q$ gives the framing of $Q$. Additionally, $Q$ is $G$-invariant and codimension 0 inside of $M$, so by \ref{prop:codim zero cobordism}, $(M,f)=(Q,q)$ in $\omega_V^G[F,F'](X,A)$.
\end{proof}

\begin{prop}\label{prop:omega pi isomorphism}
The map
\[
\omega_V^G[F,F'](X,A) \xrightarrow{\Psi_\omega^H} \omega_{V^H}^{WH}[\free](X^H,A^H)
\]
is an isomorphism. 
\end{prop}

\begin{proof}
The claim follows immediately from \cref{prop:right inverse} and \cref{prop:left inverse}.
\end{proof}

\beforesubsection
\subsubsection{The map $\Psi_\pi^H$ is an isomorphism}\aftersubsection

Now that we have shown that the map $\Psi_\omega^H$ on cobordism groups is an isomorphism, we will show that the map $\Psi_\pi^H$ on stable homotopy is an isomorphism as well. 

\begin{df}\label{df:geometric fixed points}
Let $H\leq G$ and let $X$ be a cofibrant genuine $G$-spectrum. Define the $G$-space $\widetilde{EP_H}$ up to homotopy as any $G$-CW complex with fixed points of the following type:
\[
\widetilde{EP_H}^K \simeq\begin{cases}
S^0 & \text{ if } H\leq K\\
* & \text{ else. }
\end{cases}
\]
Then the \textbf{geometric fixed points of $X$} are 
\[
(\widetilde{EP_H} \wedge X)^{H}.
\]
Here $\wedge$ denotes smashing $\widetilde{EP_H}$ with every spectrum level of $X$, and $(-)^H$ denotes the genuine fixed points of \cref{df:genuine and categorical fixed points}. 
\end{df}

Note that, by Elmendorf's theorem (\cite{elmendorf83}), the above formula uniquely determines $\widetilde{EP_H}$ up to $G$-homotopy equivalence.

\begin{lem}\label{lem:quotient space}
Let $F$ and $F'$ be two families of subgroups of $G$ which differ by the conjugacy class of $(H)$ i.e.
\[
F= F'\cup (H)
\]
There is a weak equivalence of $G$-spaces
\[
EF/EF' \simeq G\wedge_{NH} EWH_+ \wedge \widetilde{EP_H}
\]
\end{lem}

\begin{proof}
First recall that the equivariant homotopy type of $EF$ is defined by 
\[
(EF)^K = \begin{cases}
* & \text{ if } K\in F\\
\emptyset & \text{ else }
\end{cases}
\]
Thus, the equivariant homotopy type of $EF/EF'$ is given by 
\[
(EF/EF')^K = \begin{cases}
S^0 & \text{ if } K = H\\
* & \text{ else }
\end{cases} 
\]

Note that the isotropy group of any non-basepoint element of $G\wedge_{NH}EWH_+$ is conjugate to $H$. Moreover, 
\[
(G\wedge_{NH}EWH_+)^H \cong NH\wedge_{NH}EWH_+ \simeq S^0
\]
where the last map is a weak equivalence of non-equivariant spaces. In particular, we have now obtained
\[
(G\wedge_{NH}EWH_+)^K = \begin{cases}
S^0 & \text{ if } K=H\\
* & \text{ if $K$ properly contains $H$}
\end{cases}
\]
combining this with the definition of $\widetilde{EP}_H$ from \cref{df:geometric fixed points} we obtain
\[
(G\wedge_{NH}EWH_+ \wedge\widetilde{EP}_H)^K = \begin{cases}
S^0 & \text{ if } K = H\\
* & \text{ else }
\end{cases}
\]
Therefore, by Elmendorf's theorem, there is a weak equivalence of $G$-spaces
\[
G\wedge_{NH}EWH_+ \wedge\widetilde{EP}_H \simeq EF/EF'
\]
\end{proof}

\begin{lem}\label{lem:geometric fixed points of loops}
Let $(X,A)$ be a $G$-topological pair and $H$ a subgroup of $G$. Then there is a stable equivalence of $WH$-spectra
\[
\Phi^H(\Omega^V\Sigma^\infty X/A)\simeq \Omega^{V^H}\Sigma^\infty X^H/A^H
\]
\end{lem}
\begin{proof}
First note that there is a stable equivalence $\Omega^V(-) \xleftarrow{\simeq} F_VS^0\wedge(-$) and a natural isomorphism $\Phi^H(F_V (-))\cong F_{V^H}(\Phi^H(-)$. Moreover, $\Phi^H$ commutes with smash products of cofibrant spectra \cite{mandell2002equivariant}.

Therefore, 
\begin{equation}\label{eq:geometric fixed points of loops}
\Phi^H(\Omega^V\Sigma^\infty X/A) \simeq \Omega^{V^H}\Phi^H(\Sigma^\infty X/A)
\end{equation}
Since geometric fixed points commute with suspension spectra and fixed points commute with pushouts along closed inclusions, \cref{eq:geometric fixed points of loops} becomes
\[
\Phi^H(\Omega^V\Sigma^\infty X/A)\simeq \Omega^{V^H}\Sigma^\infty X^H/A^H.
\]
\end{proof}

\begin{prop}\label{prop:psi pi isomorphism}
The map 
\[
\Psi_\pi^H :\pi_V^G[F,F'](\Sigma^\infty X/A) \rightarrow \pi_{V^H}^{WH}[\free](\Sigma^\infty X^H/A^H)
\]
is an isomorphism.
\end{prop}

\begin{proof}
First observe that 
\[
\pi_V^G[F,F'](\Sigma^\infty X/A) \cong \pi_0((\Omega^V\Sigma^\infty EF/EF'\wedge X/A)^G)
\]
So we will analyze the spectrum
\[
(\Omega^V\Sigma^\infty EF/EF'\wedge X/A)^G.
\]

We begin with the observation that 
\[
\Omega^V\Sigma^\infty EF/EF'\wedge X/A \simeq EF/EF'\wedge \Omega^V\Sigma^\infty X/A
\]
and therefore, by \cref{lem:quotient space} we have
\begin{equation}\label{eq:induction equation}
\Omega^V\Sigma^\infty EF/EF'\wedge X/A \simeq G\wedge_{NH}EWH_+\wedge \widetilde{EP}_H\wedge \Omega^V\Sigma^\infty X/A.
\end{equation}
Now taking the genuine $G$-fixed points of both sides of \cref{eq:induction equation} and applying the Wirthm\"uller isomorphism \cite[3.2.15]{schwede_global} we obtain

\begin{align*}
(\Omega^V\Sigma^\infty EF/EF'\wedge X/A)^G
& \simeq (G\wedge_{NH}EWH_+\wedge \widetilde{EP}_H\wedge \Omega^V\Sigma^\infty X/A)^G\\
& \simeq (S^L\wedge EWH_+\wedge \widetilde{EP}_H\wedge\Omega^V\Sigma^\infty X/A)^{NH}\\
& \simeq ((EWH_+\wedge \widetilde{EP}_H\wedge\Omega^V\Sigma^\infty X/A)^H)^{WH}\\
& \simeq \Phi^H(EWH_+\wedge \Omega^V\Sigma^\infty X/A)^{WH}.
\end{align*}

The second equivalence is the Wirthm\"uller isomorphism \cite[3.2.15]{schwede_global}. The third is given by iterating fixed points and observing that since $G$ is a product of a finite group and a torus and since $L$ is the $NH$-representation $T_{eNH}G/NH$, it is 0-dimensional. The final identity is the definition of geometric fixed points. 

Now by \cref{lem:geometric fixed points of loops} this becomes
\[
(\Omega^V\Sigma^\infty EF/EF'\wedge X/A)^G \simeq (EWH_+\wedge\Omega^{V^H}\Sigma^\infty X^H/A^H)^{WH}.
\]
But this is exactly the statement that 
\begin{equation}\label{eq:fixed point iso}
\pi_V^G[F,F'](\Sigma^\infty X/A)\cong \pi_{V^H}^{WH}[\free](\Sigma^\infty X^H/A^H).
\end{equation}
Recall that $\Psi_\pi$ was defined by first taking geometric $H$-fixed points and then taking genuine $WH$-fixed points. This is exactly the isomorphism we constructed in this proof. Therefore, $\Psi_\pi$ is the isomorphism of \cref{eq:fixed point iso}. 
\end{proof}

We will now assemble the preceding results together to prove an equivariant Pontryagin-Thom isomorphism. We will begin by working with a pair of families before moving to the general case. 

\begin{prop}\label{prop:adjacent families isomorphism}
Let $F$ and $F'$ be adjacent families of subgroups of $G$. The map
\[
\omega_V^G[F,F'](X,A)\xrightarrow{\Delta} \pi_V^G[F,F']\left( \Sigma^\infty X/A \right)
\]
of \cref{dia:isotropy separation diagram} is an isomorphism.
\end{prop}

\begin{proof}
The result follows immediately from \cref{prop:omega pi isomorphism} \cref{prop:free PT isomorphism} and \cref{prop:psi pi isomorphism}.
\end{proof}

\begin{lem}\label{lem:subconjugation}
Let $G$ be a compact Lie group, $H\leq G$ a closed subgroup, and $g\in G$. If $gHg^{-1}\leq H$ then $gHg^{-1}=H$. 
\end{lem}

\begin{proof}
First suppose that $H$ is connected. Since conjugation is a diffeomorphism, $\Lie(H)$ and $\Lie(gHg^{-1})$ have the same dimension where $\Lie(-)$ is the associated Lie algebra of a Lie group. Since $\Lie(gHg^{-1})\leq \Lie(H)$ the Lie algebras must be equal. Therefore, the following square commutes:
\begin{diagram*}
\Lie(gHg^{-1})\arrow[r,"\cong"]\arrow[d] & \Lie(H)\arrow[d]\\
gHg^{-1}\arrow[r]&H
\end{diagram*}
where the bottom map is inclusion and the vertical maps are the exponential. Since $G$, and thus $H$, are compact, and $H$ is connected, the exponential map is a surjection and the claim follows immediately. 

Now consider the case where $H$ is not connected. Since conjugation is a diffeomorphism, $H$ and $gHg^{-1}$ have the same number of components. Moreover, each component of $gHg^{-1}$ is contained in exactly one component of $H$. By repeating the above argument in each component we obtain the desired result. Note that in the non-identity components we consider a tangent space at a point in $gHg^{-1}$ rather than the associated Lie algebra. This difference has no effect on the proof. 
\end{proof}

\begin{lem}\label{lem:families of subgroups are adjacent or colimits}
Recall that two families of subgroups of $G$ are adjacent if they differ by a single conjugacy class. Given a family of subgroups $F$ of $G$, either $F$ is adjacent to a subfamily of $F$, or $F$ is the colimit of all families properly contained in $F$. 
\end{lem}

\begin{proof}
First observe that $\underset{F_i\subsetneq F}{\colim} F_i$ is a family of subgroups of $G$. Now suppose that 
\[
F\neq \underset{F_i\subsetneq F}{\colim} F_i.
\] 
Since $\underset{F_i\subsetneq F}{\colim} F_i \subsetneq F$, it must be the case that there is some subgroup $H$ of $G$ with
\[
H\in F - \underset{F_i\subsetneq F}{\colim} F_i.
\]
We claim that $H$ is minimal in $F - \underset{F_i\subsetneq F}{\colim} F_i$. To see this, begin by defining $F'$ to be the family constructed by taking all proper subgroups of all conjugates of $H$.  By \cref{lem:subconjugation} $H$ is not in $F'$. Therefore, every proper subgroup of $H$ is in a family properly contained in $F$ so that $H$ is minimal as desired.

Since $H$ is minimal in $F-\underset{F_i\subsetneq F}{\colim} F_i$, it must be the case that $\underset{F_i\subsetneq F}{\colim} F_i\cup (H)$ is a family of subgroups of $G$. We immediately obtain 
\begin{equation}\label{eq:containment}
\underset{F_i\subsetneq F}{\colim} F_i\cup (H)\subset F
\end{equation}
if this containment were proper, then $H$ would have to be an element of $\underset{F_i\subsetneq F}{\colim} F_i$ which is not true by assumption. Therefore, the containment of \cref{eq:containment} cannot be proper so that 
\[
F = \underset{F_i\subsetneq F}{\colim} F_i\cup (H)
\] 
as desired. 
\end{proof}

\begin{thm}\label{thm:equivariant framed PT}
The map 
\[
\omega_V^G(X,A) \xrightarrow{\Delta} \pi_V^G(\Sigma^\infty X/A)
\]
of \cref{dia:big diagram}, is an isomorphism.
\end{thm}

\begin{proof}
Let $s(G)$ be the poset of families of subgroups of $G$ ordered by inclusion. We denote $F_0=\emptyset$ and $F_1=\{e\}$. Let $t(G)$ be an extension of $s(G)$ to a total order. By \cref{prop:adjacent families isomorphism}, the map
\[
\omega_V^G[F_1,F_0](X,A)\xrightarrow{\Delta} \pi_V^G[F_1,F_0](\Sigma^\infty X/A)
\]
is an isomorphism. 

Now let $F\in t(G)$ and suppose that for all $F_k<F$ 
\[
\omega_V^G[F_k,F_0](X,A)\cong \pi_V^G[F_k,F_0](\Sigma^\infty X/A)
\]

By \cref{lem:families of subgroups are adjacent or colimits} either $F$ is adjacent to a subfamily of $F$ or it is the colimit of all subfamilies properly contained in $F$. 

Assume that $F$ is adjacent to one of its subfamilies, $F'$. Then we obtain a map of long exact sequences associated to the triple $F_0<F'<F$:
\begin{figure}[H]
\center
\begin{tikzcd}[column sep = small]
\dots \arrow[r] & \omega_V^G[F',F_0](X,A) \arrow[d] \arrow[r] & \omega_V^G[F,F_0](X,A)\arrow[d]\arrow[r] & \omega_V^G[F,F'](X,A)\arrow[d]\arrow[r]&\dots\\
\dots \arrow[r] & \pi_V^G[F',F_0](\Sigma^\infty X/A) \arrow[r] & \pi_V^G[F,F_0](\Sigma^\infty X/A)\arrow[r] & \pi_V^G[F,F'](\Sigma^\infty X/A)\arrow[r]&\dots
\end{tikzcd}
\end{figure}
The left vertical map is an isomorphism by the inductive assumption and the right vertical map is an isomorphism by \cref{prop:adjacent families isomorphism}. Thus, by the five lemma, the middle vertical map is also an isomorphism.

Now suppose that instead $F= \underset{F_i\subsetneq F}{\colim}F_i$. In this case
\begin{align*}
\omega_V^G[F,F_0](X,A)
& \cong \omega_V^G[\underset{F_i\subsetneq F}{\colim}F_i,F_0](X,A)\\
& \cong \underset{F_i\subsetneq F}{\colim}\omega_V^G[F_i,F_0](X,A)\\
& \cong \underset{F_i\subsetneq F}{\colim} \pi_V^G[F_i,F_0](\Sigma^\infty X/A)\\
& \cong \pi_V^G[\underset{F_i\subsetneq F}{\colim}F_i,F_0](\Sigma^\infty X/A)\\
& \cong \pi_V^G[F,F_0](\Sigma^\infty X/A)
\end{align*}

The second isomorphism is due to the fact that bordism groups (the elements of which are compact manifolds) commute with filtered colimits. The third isomorphism is due to the inductive assumptions. The fourth isomorphism is due to the fact that equivariant stable homotopy groups, suspension spectra, and quotients all commute with filtered colimits. 

We have now shown, by induction, that for any family of subgroups, $F$, of $G$ we have
\[
\omega_V^G[F,F_0](X,A)\cong \pi_V^G[F,F_0](\Sigma^\infty X/A)
\]
In particular, taking $F= \All$, we obtain that
\[
\omega_V^G(X,A) \xrightarrow{\Delta} \pi_V^G(\Sigma^\infty X/A)
\]
is an isomorphism. 
\end{proof}

\section{The parameterized equivariant version}

Throughout this section, $G$ is a product of a finite group and a torus and $B$ is a closed $G$-manifold. In this section, we will generalize the previous results to families of $V$-framed $G$-manifolds parameterized over $B$. We will relate the cobordism groups of such families to the equivariant stable homotopy groups of a particular spectrum of sections. 

\subsection{Parameterized equivariant framed bordism}\aftersubsection

\begin{df}
Let $(X,A)$ be a $G$-topological pair as in \cref{df:pair}. Let $B$ be a closed $G$-manifold and $\psi:X\to B$ a $G$-equivariant map.  A \textbf{singular $G$-manifold over $\psi$} is a triple, $(M,f,\phi)$ where $M$ is a compact $G$-manifold, $f:(M,\partial M)\to (X,A)$ and $\phi:M\to B$ are $G$-equivariant maps such that the following diagram commutes:

\begin{figure}[H]
\center
\begin{tikzcd}
(M,\partial M)\arrow[dr, "\phi"'] \arrow[rr,"f"] & & (X,A)\arrow[dl,"\psi"]\\
& B & 
\end{tikzcd}
\end{figure}
\end{df}

\begin{rmk}
Given $(M,f,\phi)$ a singular $G$-manifold over $\psi:X\to B$ we obtain a map of retractive $G$-spaces over $B$ by taking
\begin{figure}[H]
\center
\begin{tikzcd}
(M_{+B},\partial M_{+B})\arrow[dr, "\phi"'] \arrow[rr,"f"] & & (X_{+B},A_{+B})\arrow[dl,"\psi"]\\
& B & 
\end{tikzcd}
\end{figure}
where $M_{+B}$ denoted the space $M\amalg B$. 
\end{rmk}

\begin{df}
Let $M$ be a compact $G$-manifold and $\phi:M\to B$ a $G$-equivariant map. A \textbf{$V$-framing of $M$ over $B$} is an equivalence class of a $G$-homotopy class of $G$-bundle isomorphisms
\[
TM\oplus \varepsilon_M(\R^{k})\cong \phi^*(TB)\oplus \varepsilon_M(V\oplus\R^k)
\]
where the equivalence relation is the same as in \cref{df:V-framing}.
\end{df}

\begin{df}
A \textbf{$V$-framed cobordism element over $\psi:X\to B$} is a singular $G$-manifold over $\psi$ equipped with a $V$-framing over $B$. 
\end{df}

\begin{df}
Two $V$-framed cobordism elements $(M,f,\phi)$ and $(M',f'\phi')$ over $\psi:X\to B$ are \textbf{$V$-framed cobordant over $\psi:X\to B$} if there exists, $(W,q,p)$, a singular $G$-manifold over $\psi:X\to B$ such that:
\begin{description}
\item[(i)] $W$ is a $G$-manifold with a $(V\oplus \R)$-framing over $B$.
\item[(ii)] $W$ is a $G$-cobordism between $M$ and $M'$ as in \cref{df:G-cobordism}. 
\item[(iii)] $M\cup M'\subset \partial W$ such that the induced $V$ framing on $\partial W$ restricted to $M$ and $M'$ agree with the framings on $M$ and $M'$ respectively.
\item[(iv)] $q:W\to X$ is an equivariant map such that $q|_{M}=f$ and $q|_{M'}=f'$ and $q$ maps the sides of $W$ into $A$. 
\item[(v)] $p:W\to B$ is an equivariant map restricting to $\phi$ and $\phi'$ on $M$ and $M'$ respectively.  
\end{description}
We will frequently shorten $V$-framed cobordism over $\psi:X\to B$ to cobordism when the context is clear.
\end{df}

\begin{lem}
$V$-framed $G$-cobordism over $\psi:X\to B$ is an equivalence relation. 
\end{lem}

\begin{proof}
The proof is largely the same as in \cref{lem:eq cobordism is eq rel}. The argument concerning $\phi$ proceeds the same as the argument concerning $f$.  
\end{proof}

\begin{notn}
Define $\omega_{V}^G((X,A)\xrightarrow{\psi} B)$ to be the abelian group of $V$-framed $G$-manifolds over $\psi:X\to B$ up to cobordism. Note that the elements of $\omega_{V}^G((X,A)\to B)$ are manifolds of dimension $\dim(V)+\dim(B)$. However, it is more accurate to think of them as families of manifolds of dimension $\dim(V)$ parameterized over $B$.
\end{notn}

\beforesubsection
\subsection{Parameterized spectra}\aftersubsection

In this section we briefly introduce parameterized equivariant stable homotopy theory. Further background can be found in \cite{ms} and \cite{malkiewich_convenient}.

Given a space $B$, a retractive space over $B$ is a space $X$ equipped with maps $i:B\to X$ and $p:X\to B$ such that $p\circ i$ is the identity on $B$. If $X$ is a retractive space over $B$ and $b\in B$, we denote the fiber over $b$ as $X_b$. A map between retractive spaces is a map that is compatible with the retractions in both the source and target. We denote the category of retractive space over $B$ by $\mathcal{R}(B)$. The category of retractive spaces over every possible base is denoted $\mathcal{R}$. The category $\mathcal{R}$ is symmetric monoidal with respect to the external smash product $\barsmash$ defined in \cite{malkiewich_convenient}. If $X\in \mathcal{R}(A)$ and $Y\in \mathcal{R}(B)$ then $X\barsmash Y\in \mathcal{R}(A\times B)$ and has fiber over $(a,b)$ given by $X_a\wedge Y_b$. If $A=*$ we may consider $X\barsmash A$ as an object of $\mathcal{R}(B)$ 

\begin{ex}
Let $S^1$ and $X$ be  retractive spaces over $*$ and $B$ respectively. Define the fiberwise reduced suspension of $X$ over $B$ as 
\[
\Sigma_B X = S^1\barsmash X.
\]
\end{ex}

Parameterized orthogonal spectra form a category which we denote $\Osp(B)$.

\begin{df}
Given a continuous map $f:A\to B$ we define the \textbf{pullback} to be the functor $f^*:\mathcal{R}(B)\to \mathcal{R}(A)$ which send a retractive space over $B$ to its pullback along $f$. Note that the pullback of $X\in \mathcal{R}(B)$ along $*\to B$ is the fiber of $X$ over $*$. The external smash product commutes with pullbacks. Thus, given an orthogonal spectrum over $B$, and $b\in B$ the fibers $X_b$ form an orthogonal spectrum which we call the \textbf{fiber spectrum at $b$}. 
\end{df}

\begin{df}\label{df:spectrum of sections}
Let $X\in \Osp(B)$. The spectrum of sections of $X$, denoted $\Gamma_BX\in \Osp$ is defined by $(\Gamma_BX)_n = \{s:B\to X_n\mid p_n\circ s=i_n\}$ with bonding maps given by commuting this construction with the adjunct structure maps. Here $p_n$ is the projection $X_n\to B$ and $i_n:B\to X_n$ is the zero section.   
\end{df}

The spectrum of sections defines a functor $\Gamma_B:\Osp(B)\to \Osp$. This functor has a left adjoint called the pullback. 

A retractive $G$-space over a $G$-space $B$ is defined the same as a retractive space except that we require the inclusion and projection maps be $G$-equivariant. Note however, that the fibers $X_b$ of a retractive $G$-space $X$ over $B$ are $G_b$-spaces rather than $G$-spaces. We denote the category of retractive $G$-spaces by $G\mathcal{R}(B)$. Similarly, we denote the category of parameterized orthogonal $G$-spectra over $B$ by $G\Osp(B)$ and remark that the fiber spectrum $X_b$ is an orthogonal $G_b$-spectrum.

\begin{df}
For $X\in G\Osp(B)$, the stable homotopy groups of $X$ are the genuine stable homotopy groups of the orthogonal $G_b$-spectra $X_b$ for each $b\in B$. These are defined to be
\[
\pi_{n,b}^H(X) = \underset{V\subset \mathcal{U}}{\colim}\pi_n([\Omega^VX_b(V)]^H)
\]
where $n\geq 0$, and $H$ is a closed subgroup of $G_b$, and $V$ run over the representations in the universe $\mathcal{U}$. 
\end{df}

Note that the stable homotopy groups do not preserve level-wise equivalences. For example, the map induced on suspension spectra by $\{0\}_{+I}\to I_{+I}$ over $I$ is a level equivalence but does not induce an isomorphism on $\pi_{0,1/2}(-)$. However, the stable homotopy groups are right deformable. We may therefore define the right derived homotopy groups to be 
\[
\mathbb{R}\pi_{n,b}^H(X) = \pi_{n,b}^H(R^{lv}(X))
\]
where $R^{lv}$ is defined as in \cite[4.2.5]{malkiewich_convenient}

We include the following remark to clarify the definition of right derived homotopy groups. 

\begin{rmk}
We say a functor $F:C\to D$ between categories with weak equivalences is right deformable if there is a full subcategory $A\subset C$ on which $F$ preserves weak equivalences, a `fibrant replacement' functor $R:C\to C$ the image of which lands in $A$, and a natural weak equivalence $r:X\xrightarrow{\sim} RX$. In this case we define the right derived functor of $F$ as $\R F=F\circ R$. Note that left deformable and left derived functors are defined similarly except the natural weak equivalence goes from the cofibrant replacement to the identity. Further details can be found in \cite{malkiewich_convenient} and \cite{dhks}.
\end{rmk}

\begin{prop}\cite[7.4.2]{malkiewich_convenient}\label{prop:model for parameterized G-spectra}
There is a proper model structure on $G\Osp(B)$ where the weak equivalences are the maps inducing isomorphisms on the derived stable homotopy groups. 
\end{prop}

We use this model structure to work with parameterized spectra homotopically, that is, up to weak equivalence. For example, if $B$ is a cell complex, then the spectrum of sections over $B$, denoted $\Gamma_B(-)$, is a right Quillen functor by \cite[Lemma 5.4.4]{malkiewich_convenient}. We use this model structure to construct the right derived functor of $\Gamma_B$, which we denote $\R\Gamma_B$. The right derived functor $\R\Gamma_B$ is defined by first applying fibrant replacement and then applying sections. 

\begin{rmk}
Given $X\in G\Osp(B)$ the equivariant spectrum of sections is defined the same as in \cref{df:spectrum of sections}. The categorical and genuine fixed points of $X\in G\Osp(B)$ are defined analogously to the non-parameterized case. Note however that the genuine and categorical $H$-fixed points of $X\in G\Osp(B)$ are objects in $WH\Osp(B^H)$. 
\end{rmk}

\beforesubsection
\subsection{Proof of the parameterized equivariant version}\aftersubsection

We will now prove a version of the Pontryagin-Thom theorem that relates families of $V$-framed $G$-manifolds parameterized over $B$ to parameterized equivariant stable homotopy theory. 

\beforesubsection 
\subsubsection{Isotropy separation}\aftersubsection

We now build a commutative diagram analogous to that in \cref{prop:big diagram commutes}.

\begin{df}\label{df:parameterized equivariant framed PT}
The \textbf{parameterized equivariant framed Pontryagin-Thom construction} denoted $\Delta$ is a map
\[
\omega_{V}^G((X,A)\xrightarrow{\psi} B) \xrightarrow{\Delta} \pi_V^G\left( \R\Gamma_B(\Sigma_B^\infty X\cup_A B) \right)
\]
defined as follows.

Let $(M,f,\phi)\in \omega_{V}^G((X,A)\to B)$. Since $M$ is $V$-framed over $B$, we have
 
\begin{equation}\label{eq:framing over B}
TM\oplus\varepsilon_M(\R^{k})\cong \phi^*(TB)\oplus\varepsilon_M(V\oplus\R^k)
\end{equation}

Let $\rho$ be the regular representation of $G$, and embed $(M,\partial M)\xrightarrow{i} (D(n\rho\oplus V),S(n\rho\oplus V))$. We then obtain a fiberwise embedding
\[
(M,\partial M)\xrightarrow{(\tau,i)} B\times(n\rho\oplus V)
\]
Call the normal bundle of this embedding $\nu$.

We may add $\nu$ to both sides of \cref{eq:framing over B} to obtain
\begin{align*}
\nu\oplus \tau^*(TB)\oplus \varepsilon(V\oplus\R^k)
& \cong \nu\oplus TM\oplus\varepsilon(\R^{k})\\
& \cong (\tau,i)^*(T(B\times n\rho\oplus V))\oplus\varepsilon(\R^{k})\\
& \cong \tau^*(TB)\oplus \varepsilon(n\rho\oplus V\oplus \R^{k})
\end{align*}

By adding the inverses of $\tau^*(TB)$ and $\varepsilon(V)$ and increasing the value of $k$, we obtain
\begin{equation}\label{eq:trivial normal bundle}
\nu\oplus\varepsilon(\R^k)\cong \varepsilon(n\rho\oplus\R^{k}).
\end{equation}

We fiberwise one point compactify $B\times(n\rho\oplus V)$ then apply a fiberwise Pontryagin-Thom collapse map to obtain
\[
B\times S^{n\rho\oplus V} \xrightarrow{c} D(\nu)\cup_{D(\nu|_{\partial M})\cup S(\nu)}B.
\]
Using \cref{eq:trivial normal bundle}, we have a map,
\[
D(\nu)\cup_{D(\nu|_{\partial M})\cup S(\nu)}B \xrightarrow{\eta} (M\times D(n\rho))\cup_{(\partial M\times D(n\rho))\cup (M\times S(n\rho))}B
\]
Now, applying $f:(M,\partial M)\to (X,A)$ we obtain
\[
 M\times D(n\rho)\cup_{\partial M\times D(n\rho)\cup M\times S(n\rho)}B \xrightarrow{\tilde{f}}  X\times D(n\rho)\cup_{A\times D(n\rho)\cup X\times S(n\rho)}B
\]

The composition of these three maps defines
\[
B\times S^{n\rho\oplus V} \xrightarrow{\tilde{f}\circ \eta\circ c} X\times D(n\rho)\cup_{A\times D(n\rho)\cup X\times S(n\rho)}B
\]
Finally, we can simply expand the definitions of $(X\times D(n\rho))\cup_{(A\times D(n\rho))\cup (X\times S(n\rho))}B$ and $(X\cup_A B) \barsmash S^{n\rho}$ and see that they are the same on the level of sets. It is then straightforward to check that they are the same as $G$-spaces over $B$. In particular, by identifying $(X\times D(n\rho))\cup_{(A\times D(n\rho))\cup (X\times S(n\rho))}B$ and $X\cup_A B \barsmash S^{n\rho}$ we see that $\tilde{f}\circ \eta\circ c$ defines an element of $\pi_V^G(\Gamma_B\Sigma_B^\infty X\cup_A B)$. Composing with the map $\Gamma_B(-)\to \R\Gamma_B$ gives an element of the desired group.

Note that in the above, $X\cup_A B$ is the fiberwise analogue of collapsing $A$ to the basepoint. 

We call this assignment of $(M,\partial M)$ to $\tilde{f}\circ \eta\circ c$ the Pontryagin-Thom map and denote it by $\Delta$. 
\end{df}

\begin{df}
Let $H$ be a closed subgroup of $G$ and $(M,f,\phi)\in \omega_V^G((X,A)\xrightarrow{\psi}B)$. Then $(M^H,f^H,\phi^H)$ is a $WH$-manifold which is $V^H$-framed over $\psi^H:X^H\to B^H$. Define the $H$-fixed point map as 
\begin{align*}
\omega_V^G((X,A)\xrightarrow{\psi}B)&\xrightarrow{\Psi_\omega^H} \omega_{V^H}^{WH}((X^H,A^H)\xrightarrow{\psi^H}B^H)\\
(M,f,\phi)& \mapsto (M^H,f^H,\phi^H)
\end{align*}
\end{df}

\begin{df}
Let $H$ be a closed subgroup of $G$. We define
\[
\pi_V^G \left(\Gamma_BR^{lv}( \Sigma^\infty X\cup_A B) \right) \xrightarrow{\Psi_\pi^H} \pi_{V^H}^{WH}\left( \Gamma_BR^{lv}(\Sigma_B^\infty X^H\cup_{A^H} B^H) \right)
\]
to be the genuine $H$-fixed points of the geometric $H$-fixed points. In other words, we send the $G$-homotopy class of the $G$-equivariant map 
\[
B\times S^{k\rho\oplus V}\rightarrow S^{k\rho}\wedge R^{lv}(X\cup_A B)
\]
to the $WH$-homotopy class of the $WH$-equivariant map 
\[
B^H\times S^{k\rho^H\oplus V^H}\rightarrow S^{k\rho^H}\wedge R^{lv}(X\cup_{A} B)^H
\]
where $R^{lv}$ is the levelwise fibrant replacement functor of \cite{malkiewich_convenient}. In \cref{lem:derived sections}, we will show that $\Gamma_B(R^{lv}(-))$ agrees with the right derived functor of $\Gamma_B(-)$ in the stable model structure. 
\end{df}

\begin{df}
Finally, we define a map 
\[
\omega_{V^H}^{WH}((X^H,A^H)\xrightarrow{\psi^H} B^H) \rightarrow \pi_{V^H}^{WH}(\Gamma_{B^H}R^{lv}(\Sigma^\infty_{B^H} X^H\cup_{A^H}B^H))
\]
by applying the construction in \cref{df:parameterized equivariant framed PT} to $V^H$-framed $WH$-cobordism classes over $\psi^H:X^H\to B^H$. 
\end{df}

\begin{notn}
Let $F$ and $F'$ be two families of subgroups of $G$ with $F'\subset F$. We define
\[
\omega_V^G[F,F']((X,A)\to B) = \omega_V^G((EF\times X, EF\times A\cup EF'\times X)\to B)
\]
and
\[
\pi_V^G[F,F'] (\Gamma_BR^{lv}(\Sigma^\infty_B X\cup_A B)) = \pi_V^G(\Gamma_BR^{lv}(\Sigma^\infty_B EF\times X\cup_{EF'\times X \cup EF\times A} B)).
\]
Note that $EF\times X\cup_{EF'\times X \cup EF\times A} B\cong EF/EF'\barsmash X\cup_A B$.
\end{notn}

\begin{prop}
Let $H$ be a closed subgroup of $G$ and let $F$ and $F'$ be two families of subgroups of $G$ such that $F=F'\cup (H)$. Then the following diagram commutes:

\begin{diagram}\label{dia:parameterized big diagram}
\omega_{V}^G[F,F']((X,A)\xrightarrow{\psi} B) \arrow[r,"\Delta"] \arrow[d,"\Psi_\omega^H"] & \pi_V^{G}[F,F'](\Gamma_BR^{lv}(\Sigma^\infty_B X\cup_A B))\arrow[d,"\Psi_\pi^H"]\\
\omega_{V^H}^{WH}[\free]((X^H,A^H)\xrightarrow{\psi^H}B^H) \arrow[r,"\Delta^H"] & \pi_{V^H}^{WH}[\free](\Gamma_{B^H}R^{lv}(\Sigma^\infty_{B^H} X^H\cup_{A^H} B^H)).
\end{diagram}
\end{prop}

\begin{proof}
As in \cref{prop:big diagram commutes} the proof is a straightforward verification. 
\end{proof}

\beforesubsection
\subsubsection{The map $\Psi_\omega^H$ is an isomorphism}\aftersubsection

\begin{df}\label{df:manifold inverse}
Define a map 
\[
\Theta: \omega_{V^H}^{WH} ((X^H,A^H)\to B^H) \to \omega_{V}^G(X,A)
\]
by sending $(N,t,\tau)$ to $(Q,f,\phi)$ where 
\[
Q=G\times_{NH}(N\times D(V^H)^\perp\times_{B^H}D(\nu(i)))
\]
where $i:B^H\to B$ is the inclusion and $\nu(i)$ is the associated normal bundle. The map $f:Q\to X$ is defined as in \cref{df:the manifold Q}.

Given a map $\tau:N\to B^H$ we define a map $N\times D(V^H)^\perp\to B^H$ by projecting to $N$ then applying $\tau$. The map $\phi$ is then given by 
\begin{equation}\label{eq:phi tilde}
G\times_{NH}(N\times D(V^H)^\perp)\times_{B^H}D(\nu(i)) \to B^H\times_{B^H}B \simeq B.
\end{equation}
\end{df}

\begin{lem}
The $G$-manifold $Q$ of \cref{df:manifold inverse} is $V$-framed over $B$.
\end{lem}

\begin{proof}
We compute that
\begin{align*}
TQ\oplus\varepsilon_Q(\R^{k})
& \cong (TN\times \R^k) \times (D(V^{H^\perp})\times V^{H^\perp})\times_{TB^H}TB|_{D(\nu(i))}\\
& \cong \tau^*(TB^H)\times (V^H\oplus \R^k)\times (D(V^{H^\perp})\times V^{H^\perp})\times_{TB^H}TB|_{D(\nu(i))}\\
\end{align*}
first by definition and then using the $V^H$-framing of $N$ over $B^H$. 

Now, rearranging the terms and observing that one of our pullbacks is along the identity map, we obtain
\begin{align*}
TQ\oplus\varepsilon_Q(\R^{k})
&\cong N\times_{B^H}TB^H\times_{TB^H}TB|_{D(\nu(i))}\times (V^H\oplus \R^k)\times (D(V^H)^\perp\times (V^H)^\perp)\\
& \cong N\times_{B^H}TB|_{D(\nu(i))}\times (D(V^H)^\perp\times(V\oplus \R^k))
\end{align*}

Rewriting the restriction bundle as a pullback along an inclusion, we finally see that
\begin{align*}
TQ\oplus\varepsilon_Q(\R^{k})
& \cong N\times_{B^H}D(\nu(i))\times_B TB\times (D(V^H)^\perp\times(V\oplus \R^k))\\
& \cong \phi^*(TB)\oplus \varepsilon(V\oplus \R^k).
\end{align*}
\end{proof}
 
\begin{prop}\label{prop:psi theta is identity}
The composition $\Psi_\omega^H\Theta$ is the identity.
\end{prop}

\begin{proof}
The proof that $\Psi_\omega^H \Theta=1$ is the same as in \cref{prop:right inverse} with one exception. If $(N,t,\tau)\in\omega_{V^H}^{WH}((X^H,A^H)\to B^H)$, we simply need to check that $\Theta(\tau)^H=\tau$. However, this is almost immediate just by passing \cref{eq:phi tilde} to the $H$-fixed points. 
\end{proof}

\begin{prop}\label{prop:theta psi is identity}
The composition $\Theta\Psi_\omega^H$ is the identity.
\end{prop}

\begin{proof}
Let $(M,f,\phi)\in \omega_V^G((X,A)\to B)$. The proof that $\Theta\Psi_\omega=1$ is also the same as in \cref{prop:left inverse} except in the case of the map $\phi:M\to B$. We have that
\[
\phi|_{G\times_{NH}(M^H\times D(V^H)^\perp)\times_{B^H}D(\nu(i))} = G\times_{NH}\phi(-)^H.
\]
We can then use the compatibility of $\phi$ and $G\times_{NH}\phi(-)$ to build a map from the cobordism $W$ of \cref{prop:left inverse} to $B$ which restricts to these maps on the boundaries. This shows that $M$ is cobordant to
\[
G\times_{NH}(M^H\times D(V^H)^\perp) \in \omega_{V}^G[F,F']((X,A)\to B).
\]
\end{proof}

\begin{cor}\label{cor:eq param geom iso}
The map $\Psi_\omega^H$ is an isomorphism.
\end{cor}
\begin{proof}
This follows from \cref{prop:theta psi is identity} and \cref{prop:psi theta is identity} immediately. 
\end{proof}

\begin{prop}\label{prop:eq param fp PT iso}
The map $\Delta^H$ is an isomorphism.
\end{prop}
\begin{proof}
The proof of this fact is largely the same as the classical proof (see \cite{milnor_diff_top}) since the elements of $\omega_{V^H}^{WH}[\free]((X^H,A^H)\to B^H)$ are free $WH$-manifolds for which we can approximate any map, up to $WH$-homotopy, by an equivariant map transverse to a submanifold of our choosing. We will highlight the differences between the parameterized case and the classical proof here. 

In order to prove surjectivity, we begin with an element
\[
(B^H\times S^{n\rho\oplus V^H}\xrightarrow{f} S^{n\rho}\barsmash EWH_+\barsmash R^{lv}(X\cup_{A} B)^H)\in\pi_{V^H}^{WH}[\free]\left( \Gamma_{B^H} R^{lv}\Sigma^\infty_{B^H} X^H\cup_{A^H}B^H\right).
\] 
Note that by \cite{malkiewich_convenient} we may construct $R^{lv}$ to be strong symmetric monoidal so that we may smash with $EWH_+$ outside of $R^{lv}$ without changing the homotopy class of $f$. 

We post-compose with the map $EWH_+\barsmash R^{lv}(X\cup_{A}B)^H\to B^H\to *$ to obtain a an equivariant map of free $WH$-spaces,
\[
B^H\times S^{n\rho\oplus V^H}\to S^{n\rho}\barsmash X^H\cup_{A^H} B^H \to S^{n\rho}\times B \to S^{n\rho}.
\]
We then send $f$ to $f^{-1}(0)$ after equivariantly deforming $f$ to be smooth and have $0$ as a regular value. The framing on $f^{-1}(0)$ is given by choosing a positively oriented basis of $T_0 S^{n\rho}$ and pulling back along $f$ to obtain a trivialization
\[
\nu(f^{-1}(0),B^H\times S^{n\rho\oplus V^H})\oplus \varepsilon(\R^k)\cong \varepsilon(n\rho\oplus \R^k).
\]
Reversing the manipulation in \cref{eq:trivial normal bundle} gives the desired framing of $f^{-1}(0)$. 

The arguments for the rest of the proof are the same as in \cite{milnor_diff_top}.
\end{proof}

The only remaining step of the proof that
\[
\omega^G_{V}((X,A)\to B) \xrightarrow{\Delta} \pi_V^G(\Gamma_BR^{lv}\Sigma^\infty X\cup_A B)
\]
is an isomorphism, is proving that the map $\Psi_\pi^H$ of \cref{dia:parameterized big diagram} is an isomorphism. The proof of this fact relies on several long technical results which, for the sake of the reader, have been relegated to \cref{section:appendix}.

\beforesubsection
\subsubsection{The map $\Psi_\pi^H$ is an isomorphism}\aftersubsection

We now prove that the map $\Psi_\pi^H$ of \cref{dia:parameterized big diagram} is an isomorphism.

\begin{prop}\label{prop:parameterized fixed point iso}
The map 
\[
\Psi_\pi^H:\pi_V^G[F,F'](\Gamma_BR^{lv}(\Sigma^\infty_B X\cup_A B))\to \pi_{V^H}^{WH}[\free](\Gamma_BR^{lv}(\Sigma^\infty_{B^H} X^H\cup_{A^H} B^H))
\]
is an isomorphism. 
\end{prop}

\begin{proof}
First observe that 
\begin{align*}
\pi_V^G[F,F'](\Gamma_BR^{lv}(\Sigma_B^\infty X\cup_AB)) \cong \pi_0((\Omega^V\Gamma_BR^{lv}(\Sigma_B^\infty EF/EF'\barsmash X\cup_AB))^G)
\end{align*}
So we will analyze the spectrum
\[
(\Omega^V\Gamma_BR^{lv}(\Sigma_B^\infty EF/EF'\barsmash X\cup_A B))^G.
\]

By \cref{prop:sections commutes with smash} and the fact that $R^{lv}$ may be constructed to be strong monoidal as in \cite{malkiewich_convenient}, 
\[
\Omega^V\Gamma_BR^{lv}(\Sigma_B^\infty EF/EF'\barsmash X\cup_A B) \simeq EF/EF'\wedge \Omega^V\Gamma_BR^{lv}(\Sigma_B^\infty X\cup_A B).
\]
Therefore, by \cref{lem:quotient space} we have
\begin{equation}\label{eq:induction equation}
\Omega^V\Gamma_BR^{lv}(\Sigma_B^\infty EF/EF'\barsmash X\cup_A B) \simeq G\wedge_{NH}(EWH_+\wedge \widetilde{EP}_H\wedge \Omega^V\Gamma_BR^{lv}(\Sigma_B^\infty X\cup_A B)).
\end{equation}

Now taking the genuine $G$-fixed points of both sides of \cref{eq:induction equation} and applying the Wirthm\"uller isomorphism \cite[3.2.15]{schwede_global} we obtain

\begin{align*}
(\Omega^V\Gamma_BR^{lv}(\Sigma_B^\infty EF/EF'\barsmash X\cup_A B))^G
& \simeq (G\wedge_{NH}EWH_+\wedge \widetilde{EP}_H\wedge \Omega^V\Gamma_BR^{lv}(\Sigma_B^\infty X\cup_A B))^G\\
& \simeq (S^L\wedge EWH_+\wedge \widetilde{EP}_H\wedge\Omega^V\Gamma_BR^{lv}(\Sigma_B^\infty X\cup_A B))^{NH}
\end{align*}

where $L$ is the $NH$-representation $T_{eNH}(G/NH)$. 

Recall however, that since $G$ is a product of a finite group and a torus, $G/NH$ is zero dimensional so that $S^L\cong S^0$. 

Iterating fixed points and applying the definition of geometric fixed points, the above equation becomes

\[
(\Omega^V\Gamma_BR^{lv}(\Sigma_B^\infty EF/EF'\barsmash X\cup_A B))^G \simeq (\Phi^H(EWH_+\wedge\Omega^V\Gamma_BR^{lv}(\Sigma_B^\infty X\cup_A B)))^{WH}
\]

Now by \cref{prop:sections commute with geom fp} and \cref{lem:geometric fixed points of loops} this becomes
\[
(\Omega^V\Gamma_BR^{lv}(\Sigma_B^\infty EF/EF'\barsmash X\cup_A B))^G \simeq (EWH_+\wedge\Omega^{V^H}\Gamma_{B^H}R^{lv}(\Sigma_{B^H}^\infty X^H\cup_{A^H} B^H))^{WH}.
\]
But this is exactly the statement that 
\begin{equation}\label{eq:fixed point iso}
\pi_V^G[F,F'](\Gamma_BR^{lv}(\Sigma_B^\infty X\cup_A B))\cong \pi_{V^H}^{WH}[\free](\Gamma_{B^H}R^{lv}(\Sigma_{B^H}^\infty X^H\cup_{A^H} B^H))
\end{equation}
after commuting $EHW_+$ back past $\Gamma_{BH}$ and $R^{lv}$. Recall that $\Psi_\pi^H$ was defined by first taking geometric $H$-fixed points and then taking genuine $WH$-fixed points. This is exactly the isomorphism we constructed in this proof. Therefore, $\Psi_\pi^H$ is the isomorphism of \cref{eq:fixed point iso}. 
\end{proof}

We are now ready to prove the main result of the paper, the parameterized equivariant Pontryagin-Thom theorem for framed manifolds.

\begin{thm}\label{thm:equivariant parameterized PT}
Let $G$ be isomorphic to the product of a finite group and a torus. The map
\[
\omega_{V}^G((X,A)\xrightarrow{\psi} B) \xrightarrow{\Delta} \pi_V^G\left( \Gamma_BR^{lv}(\Sigma_B^\infty X\cup_A B) \right)
\]
of \cref{df:parameterized equivariant framed PT} is an isomorphism.
\end{thm}
\begin{proof}
The map 
\[
\omega_{V}^G[F,F']((X,A)\xrightarrow{\psi} B) \xrightarrow{\Delta} \pi_V^G[F,F']\left( \Gamma_BR^{lv}(\Sigma_B^\infty X\cup_A B) \right)
\]
is an isomorphism by \cref{cor:eq param geom iso}, \cref{prop:eq param fp PT iso}, and \cref{prop:parameterized fixed point iso}. The argument that 
\[
\omega_{V}^G((X,A)\xrightarrow{\psi} B) \xrightarrow{\Delta} \pi_V^G\left(\Gamma_BR^{lv}(\Sigma_B^\infty X\cup_A B) \right)
\]
is an isomorphism is then identical to the proof of \cref{thm:equivariant framed PT}.
\end{proof}

\appendix
\section{Technical results on parameterized $G$-spectra}\label{section:appendix}

This section contains a handful of technical results on parameterized $G$-spectra which are necessary for the proof of \cref{prop:parameterized fixed point iso}.

The derived functors referred to in this section are constructed from the model structure of \cref{prop:model for parameterized G-spectra}. It is however possible to derive some functors by using less than the entire (co)fibrant replacement functor of the stable model structure. We will point out when this occurs. We begin by defining several different notions of (co)fibrancy that will arise throughout this section. 

\begin{df} A map $X\to Y$ of retractive spaces over $B$ is:
\begin{enumerate}
\item an \textbf{$h$-fibration} (or Hurewicz fibration) if it has the homotopy lifting property. This means we ask for the map $X^I\to X\times_Y Y^I$ to have a section. 
\item a \textbf{$q$-fibration} (or Serre fibration) if it has the homotopy lifting property with respect to cylinders on discs:
\begin{diagram*}
D^n\times\{0\} \arrow[r] \arrow[d]& Y\arrow[d]\\
D^n\times I \arrow[r]\arrow[ru, dashed]& X
\end{diagram*}
\item an \textbf{$h$-cofibration} if it is a closed inclusion and has the homotopy extension property. In essence, the inclusion $X\times I\cup_{X\times\{0\}} Y\times\{0\}\to Y\times I$ has a retract. 
\item an \textbf{$f$-cofibration} if it a closed inclusion and has the fiberwise homotopy extension property. In essence, the inclusion $X\times I\cup_{X\times\{0\}} Y\times\{0\}\to Y\times I$ has a retract which respects the projection map to $B$. 
\end{enumerate}
\end{df}

We say a retractive space is $h$-fibrant if the map to the zero object is an $h$-fibration. A retractive space is $h$-cofibrant if the map from the zero object is an $h$-cofibration. Then (co)fibrancy for the other notions above is defined similarly. We say that a parameterized orthogonal spectrum is level $h$-fibrant if each of its constituant spaces is level $h$-fibrant. The other notions of level (co)fibrancy are defined similarly. 

\begin{df}
A parameterized orthogonal spectrum is stably fibrant if it is level $q$-fibrant and the maps $X_n \to \Omega_B X_{n+1}$ are weak equivalences for all $n\geq 0$. 
\end{df}

When defining the weak equivalences and $q$-fibrations of $G$-equivariant retractive spaces, we ask that $f:X\to Y$ be a weak equivalence or $q$-fibration on the fixed points $f^H:X^H\to Y^H$ for all closed subgroups $H\leq G$. To define $h$-fibrations, $h$-cofibrations, and $f$-cofibrations of equivariant retractive spaces, we simply require that the relevant section or retract be equivariant as well. To define level (co)fibrations of $G$-equivariant parameterized spectra, the group of equivariance at level $n$ is $G\times O(n)$ rather than just $G$. We finally note that a $G$-spectrum over $B$ is fibrant if it is level $q$-fibrant and the adjunct bonding maps are weak equivalences on all of the fixed point subspaces. 

\begin{df}
We define the \textbf{monoidal geometric fixed points} $\Phi^H_M(X)$ of a $G$-spectrum $X$ over a $G$-space $B$ to be the coequalizer  
\begin{equation}
\bigvee_{V,W} F_{W^H}S^0 \barsmash \sJ_G^H(V,W) \barsmash X(V)^H \rightrightarrows \bigvee_V F_{V^H}S^0 \barsmash X(V)^H \rightarrow \Phi_M^H(X)
\end{equation}

The action $F_{W^H}S^0 \barsmash J_G^H(V,W)\to F_{V^H}(S^0)$ is given at spectrum level $U$ by
\begin{align*}
F_{W^H}S^0(U) \barsmash J_G^H(V,W)(U) \simeq J_G(W^H,U)\barsmash J_G^H(V,W)  & \to J_G(V^H,Y)\\
(W^H\xrightarrow{i} U,V\xrightarrow{j} W)&\mapsto (V^H\xrightarrow{j^H} W^H\xrightarrow{i} Y)
\end{align*}
The other map in the coequalizer diagram is given by first considering $X$ as a $\GTop_*$ valued diagram on $J_G$ and applying it to $J_G^H(V,W)$. 
\end{df}

The monoidal geometric fixed points are a functor $\Phi_M^H(-):G\cS(B)\to WH\cS(B^H)$. 

Given a closed subgroup $H$ of $G$, recall the definition of the space $\widetilde{EP_H}$ from \cref{df:geometric fixed points}. Considering this space as a $G$ space over a point, we have the following proposition.

\begin{prop}\label{prop:geometric fixed points}
Let $X$ be a cofibrant orthogonal $G$-spectrum over a $G$-space $B$, then the geometric fixed points, $\Phi^H(X)=(f(\widetilde{EP_H}\barsmash X))^H$, agree with the monoidal geometric fixed points of $X$ up to stable equivalence, where $f(-)$ is the fibrant replacement functor from the model structure of \cref{prop:model for parameterized G-spectra}. 
\end{prop}

\begin{proof}
The following zig-zag of maps in $WH\cS(B^H)$ is a stable equivalence:
\[
(f(\widetilde{EP_H}\barsmash X))^H\xrightarrow{R}\Phi_M^H(f(\widetilde{EP_H}\barsmash X)) \xleftarrow{\alpha} \Phi_M^H(\widetilde{EP_H}\barsmash X)\xrightarrow{\beta}\Phi_M^H(X)
\]
The map $\alpha$ is given by applying geometric fixed points to the acyclic cofibration $\widetilde{EP_H}\barsmash X\to f(\widetilde{EP_H}\barsmash X)$. Note that $\Phi_M^H$ preserves acyclic cofibrations by \cite{mp2}, so $\alpha$ is also an acyclic cofibration and in particular, a weak equivalence. 

The map $\beta$ is constructed from the chain of equivalences as follows:
\begin{align*}
\Phi_M^H(\widetilde{EP_H}\barsmash X)
& \simeq \widetilde{EP_H}^H\barsmash \Phi_M^H(X)\\
& \simeq \Phi_M^H(X)
\end{align*}
and is thus a weak equivalence.

The restriction $R$ gives a map of $WH$-spectra over $B^H$. For each $b\in B^H$ the map induced by $R$ on the fiber over $B$ is the non-parameterized restriction map of $WH_b$-spectra which is an equivalence due to our cofibrancy assumption. Since $R$ gives an equivalence in each fiber it is a level equivalence. 
\end{proof}

As noted in \cite{mp2} and \cite{mandell2002equivariant} the geometric fixed points functor commutes with smash products of cofibrant spectra. Moreover, we have the following equivalence
\[
\Phi^H F_VX \simeq F_{V^H}X^H.
\]

\begin{lem}\label{lem:serre fibration}
Let $B$ be a finite $G$-CW-complex and let $r:B\to *$ be the map from $B$ to the terminal $G$-space. Let $X\in G\Osp(B)$ be an orthogonal $G$-spectrum parameterized over $B$ which is level $q$-fibrant. Define a map 
\[
(r|_{B^{(n)}})_*(X) \xrightarrow{j} (r|_{B^{(n-1)}})_*(X)
\]
by restricting a section of $X$ on $B^{(n)}$ to a section on $B^{(n-1)}$. The map $j$ is a level-wise Serre fibration. 
\end{lem}

\begin{proof}

Let $h:D^k\times [0,1]\to (r|_{B^{(n-1)}})_*X$ a homotopy and $\tilde{h_0}$ a lift of $h$ to $D^k\times \{0\} \to (r|_{B^{(n-1)}})_*X$. The following diagram commutes:
\begin{figure}[H]
\center
\begin{tikzcd}
D^k\times\{0\}\arrow[d]\arrow[r] & (r|_{B^{(n)}})_*X\arrow[d,"j"]\\
D^k\times I\arrow[r] & (r|_{B^{(n-1)}})_*X
\end{tikzcd}
\end{figure}

Rewriting $h$ and $\tilde{h_0}$ and using the commutativity of the diagram, gives us, for each $V$, the data of a map
\[
B^{(n)}\times D^k\times\{0\}\bigcup_{B^{(n-1)}\times D^k\times I} B^{(n-1)}\times D^k\times I \to X(V)
\]
This is a map of retractive spaces over $B$. There is also a map 
\[
B^{(n)}\times D^k\times\{0\}\bigcup_{B^{(n-1)}\times D^k\times I} B^{(n-1)}\times D^k\times I\to B^{(n)}\times D^k\times I
\]
over $B$ given by taking the pushout of the acyclic cofibration
\[
B^{(n)}\times D^k\times \{0\} \to B^{(n)}\times D^k\times I
\]
and the cofibration 
\[
B^{(n-1)}\times D^k\times I\to B^{(n)}\times D^k\times I
\]
Note then that this gives an acyclic cofibration over $B$. This data then assembles into a commuting square
\begin{figure}[H]
\center
\begin{tikzcd}
B^{(n)}\times D^k\times\{0\}\bigcup_{B^{(n-1)}\times D^k\times I} B^{(n-1)}\times D^k\times I \arrow[r]\arrow[d] & X(V)\arrow[d]\\
B^{(n)}\times D^k\times I \arrow[r]& B
\end{tikzcd}
\end{figure}

As previously noted, the map on the left is the pushout of an acyclic cofibration and a cofibration and is thus an acyclic cofibration. After fibrantly replacing $X$, the map on the right automatically becomes a fibration as it is the map from $X(V)$ to the terminal object. Thus, there exists a lift

\begin{figure}[H]
\center
\begin{tikzcd}
B^{(n)}\times D^k\times\{0\}\bigcup_{B^{(n-1)}\times D^k\times I} B^{(n-1)}\times D^k\times I \arrow[r]\arrow[d] & X(V)\arrow[d]\\
B^{(n)}\times D^k\times I \arrow[ru,dashed] \arrow[r]& B
\end{tikzcd}
\end{figure}

This lift then reassembles into the lift

\begin{figure}[H]
\center
\begin{tikzcd}
D^k\times\{0\}\arrow[d]\arrow[r] & (r|_{B^{(n)}})_*X\arrow[d,"j"]\\
D^k\times I\arrow[r]\arrow[ru,dashed] & (r|_{B^{(n-1)}})_*X
\end{tikzcd}
\end{figure}
as desired.
\end{proof}

\begin{lem}\label{lem:derived sections}
Let $X,Y\in G\Osp(B)$ be such that the retraction onto $B$ at each level of $X$ and $Y$ is a Serre fibration. If $X$ and $Y$ are stably equivalent, then so are $\Gamma_B X$ and $\Gamma_B Y$. 
\end{lem}

\begin{proof}
We begin by observing that the following square commutes

\begin{diagram}\label{dia:sections of serre fibrations}
\Gamma_B X \arrow[r]\arrow[d]& \Gamma_B Y\arrow[d]\\
\Gamma_B f(X) \arrow[r]& \Gamma_B f(Y)
\end{diagram}
where $f(-)$ is the stable fibrant replacement functor. 

We first claim that the bottom horizontal map is a level equivalence. Since $X\to Y$ is a stable equivalence, $f(X)\to f(Y)$ is a level equivalence. Since $B$ is a compact manifold, $\Gamma_B$ preserves level equivalences. We prove this by replacing $B$ with a finite CW-complex and inducting on the skeleta. The base case, where $B$ is a 0-dimensional CW-complex, is trivial. The inductive step is argued as follows. Assume that $\Gamma_{B^{(k)}} X(V)\simeq \Gamma_{B^{(k)}} Y(V)$. Then 
\begin{figure}[H]
\center
\begin{tikzcd}
\Gamma_{B^{(k)}\cup_{S^k}D^{k+1}} X(V) \arrow[r]\arrow[d]& \Gamma_{B^{(k)}} X(V)\arrow[d]\\
\Gamma_{D^{k+1}} X(V) \arrow[r]& \Gamma_{S^k} X(V)
\end{tikzcd}
\end{figure}
is a pullback square. Since all of the maps in the square are Serre fibrations, it is in fact a homotopy pullback square. The claim now follows from the five-lemma and the inductive assumption. 

We now claim that the vertical maps in \cref{dia:sections of serre fibrations} are stable equivalences. Observe that at level $V$, 
\begin{align*}
\Gamma_Bf(X)(V)
& \simeq \Gamma_B \hocolim_{U} \Omega_B^UX(V\oplus U)\\
& \simeq \hocolim_U \Gamma_B \Omega_B^U X(V\oplus U)\\
& \simeq \hocolim_U \Omega^U \Gamma_B X(V\oplus U)\\
& \simeq f\Gamma_B X(V).
\end{align*}
The first and last equivalences are by definition. The second equivalence is using the fact that $B$ is compact so we can commute sections over $B$ with the mapping telescope. The third equivalence follows purely from manipulating adjunctions. 

Since the vertical maps and the bottom horizontal map of \cref{dia:sections of serre fibrations} are equivalences, so is the top map. 
\end{proof}

\begin{lem}\label{lem:derived1}
Let $A$ be $h$-cofibrant and $X$ be level $h$-cofibrant and level $h$-fibrant. Then there is a zig-zag of stable equivalences 
\[
\R\Gamma_B(A\barsmash^{\mathbb{L}} X) \rightarrow \R\Gamma_B(A\barsmash X)\leftarrow \Gamma_B (A\barsmash X)
\]
\end{lem}

\begin{proof}
Since $A$ is $h$-cofibrant and $X$ is level $h$-cofibrant, the canonical map $A\barsmash^{\mathbb{L}}X \to A\barsmash X$ is a stable equivalence by \cite[Lemma 5.1.4]{malkiewich_convenient}. Moreover, $A\barsmash X$ is level $h$-fibrant by \cite[Theorem 4.4.7]{malkiewich_convenient}. Therefore, the canonical map $\Gamma_B(A\barsmash X)\to \R\Gamma_B(A\barsmash X)$ is a stable equivalence by \cref{lem:derived sections}.
\end{proof}

\begin{lem}\label{lem:derived2}
Let $A$ be $h$-cofibrant and $X$ be level $h$-cofibrant and level $h$-fibrant. Then there is a zig-zag of stable equivalences
\[
A\wedge^{\mathbb{L}} \R\Gamma_B X \leftarrow A\wedge^{\mathbb{L}}\Gamma_B X\rightarrow A\wedge \Gamma_B X. 
\]
\end{lem}

\begin{proof}
Since $X$ is level $h$-fibrant, the canonical map $\R\Gamma_B X\leftarrow \Gamma_B X$ is a stable equivalence by \cref{lem:derived sections}. Since $X$ is level $h$-fibrant and level $h$-cofibrant, it is level $f$-cofibrant by \cite{heath_kamps}. Thus, $\Gamma_B X$ is level $f$-cofibrant by \cite[Proposition 2.4.6]{malkiewich_convenient}. Since $A$ is $h$-cofibrant and $\Gamma_B X$ is level $f$-cofibrant, and thus level $h$-cofibrant, the canonical map $A\barsmash^{\mathbb{L}}\Gamma_B X\rightarrow A\barsmash \Gamma_B X$ is a stable equivalence. 
\end{proof}

\begin{lem}\label{lem:level replacement}
Every $X\in G\Osp(B)$ is stably equivalent to a parameterized $G$-spectrum over $B$ which is level $h$-cofibrant and level $h$-fibrant. 
\end{lem}

\begin{proof}
The proof follows from \cite[Proposition 4.2.5]{malkiewich_convenient}.
\end{proof}

\begin{prop}\label{prop:sections commutes with smash}
Let $B$ be a finite $G$-CW-complex. Let $X\in G\Osp(B)$ be a $G$-spectrum parameterized over $B$ and $A$ a $G$-space. Then, 
\[
A \wedge^{\mathbb{L}} \R\Gamma_B(X) \simeq \R\Gamma_B(A \barsmash^{\mathbb{L}} X)
\]
\end{prop}
\begin{proof}
By \cref{lem:derived1}, \cref{lem:derived2}, and \cref{lem:level replacement}, it suffices to assume that $A$ is $h$-cofibrant, $X$ is level $h$-cofibrant level $h$-fibrant, and then show that 
\[
A\wedge \Gamma_B X \simeq \Gamma_B(A\barsmash X).
\]

We proceed via induction on the skeleta of $B$. We first show that the result holds when restricting $r$ to the 0-skeleton of $B$. Say $B^{(0)} = \coprod G/H$ where the $H$'s may vary across closed subgroups of $G$. A section over a single orbit space $G/H$ is equivalent to the product over the elements of $G/H$ of the fiber spectra over those elements. So,
\begin{align*}
A\wedge \Gamma_{B^{(0)}}(X)
&\simeq A \wedge \prod_{b\in B^{(0)}}\Gamma_b(X)\\
&\simeq A \wedge \bigvee_{b\in B^{(0)}}\Gamma_b(X)\\
&\simeq \bigvee_{b\in B^{(0)}}A \wedge \Gamma_b(X)\\
&\simeq \prod_{b\in B^{(0)}} A \wedge \Gamma_b(X).
\end{align*}

Similarly, 
\[
\Gamma_{B^{(0)}}(A \barsmash X) \simeq \prod_{b\in B^{(0)}} \Gamma_b(A \barsmash X).
\]
Since the sections of $A \barsmash X$ are precisely the sections of $X$ smashed with $A$, we see that
\[
A\wedge \Gamma_{B^{(0)}}(X)\simeq \Gamma_{B^{(0)}}(A \barsmash X)
\]
as desired. 

We now proceed to the inductive step. Suppose that
\[
A\wedge \Gamma_{B^{(n-1)}}(X) \simeq \Gamma_{B^{(n-1)}}(A\barsmash X).
\]
By restricting a section on $B^{(n)}$ to $B^{(n-1)}$ we define a map
\[
\Gamma_{B^{(n)}}(X) \xrightarrow{j} \Gamma_{B^{(n-1)}}(X)
\]
which is a level-wise Serre fibration by \cref{lem:serre fibration}. 

The fiber $Fj$ is given level-wise by
\begin{align*}
(Fj)(V)
&=\left\{ s\in \Gamma_{B^{(n)}}(X)(V)\mid s|_{B^{(n-1)}} = i_{X(V)} \right\}\\
& = \left\{s:B^{(n)}\to X(V)\mid \pi_{X(V)}\circ s = 1 \text{ and } s|_{B^{(n-1)}} = i_{X(V)}\right\}
\end{align*}
where $i_{X(V)}$ is the zero section of $X(V)$ on $B$. 

We therefore have a fiber sequence of $G$-spectra
\[
Fj \rightarrow\Gamma_{B^{(n)}}(X) \xrightarrow{j} \Gamma_{B^{(n-1)}}(X).
\]
Smashing with $A$ preserves this fiber sequence so that we have a fiber sequence
\[
A\wedge Fj \rightarrow A \wedge\Gamma_{B^{(n)}}(X) \xrightarrow{j} A\wedge\Gamma_{B^{(n-1)}}(X).
\]
Again by \cref{lem:serre fibration}, we have a levelwise Serre fibration
\[
\Gamma_{B^{(n)}}(A\barsmash X)\xrightarrow{l} \Gamma_{B^{(n-1)}}(A\barsmash X)
\]
with fiber 
\begin{align*}
(Fl)(V) 
&=\left\{ s\in \Gamma_{B^{(n)}}(A\barsmash X)(V)\mid s|_{B^{(n-1)}} = i_{(A\barsmash X)(V)} \right\}\\
& = \left\{s:B^{(n)}\to (A\barsmash X)(V)\mid \pi_{(A\barsmash X)(V)}\circ s = 1 \text{ and } s|_{B^{(n-1)}} = i_{(A\barsmash X)(V)}\right\}.
\end{align*}
Hence, 
\[
Fl\rightarrow \Gamma_{B^{(n)}}(A\barsmash X)\xrightarrow{l} \Gamma_{B^{(n-1)}}(A\barsmash X)
\]
is a fiber sequence.

We then have a map of fiber sequences
\begin{figure}[H]
\center
\begin{tikzcd}
A\wedge Fj \arrow[d]\arrow[r] & A \wedge\Gamma_{B^{(n)}}(X) \arrow[d]\arrow[r,"j"] & A\wedge\Gamma_{B^{(n-1)}}(X)\arrow[d]\\
Fl\arrow[r] & \Gamma_{B^{(n)}}(A\barsmash X)\arrow[r,"l"] & \Gamma_{B^{(n-1)}}(A\barsmash X)
\end{tikzcd}
\end{figure}
where the two vertical maps on the right send $a\wedge (s:B^{(m)}\to X(V))$ to $\tilde{s}:B^{(m)}\to A\barsmash X(V)$ given by $\tilde{s}(b)=a\barsmash s(b)$, where $m=n,n-1$. The map on the left is the map induced on fibers. This is the map which at level $V$ is given by
\[
a\wedge s\mapsto (\tilde{s}:b\mapsto a\barsmash s(b)).
\]

By assumption, the vertical map on the right is a stable equivalence. 

We now claim that the vertical map on the left, is also a stable equivalence. Recall that the fiber $Fj$ has level $V$ given by, sections of $X(V)$ over $B^{(n)}$ which vanish on $B^{(n-1)}$. This spectrum then splits into a finite product over the $n$-cells of $B$. So $Fj$ is stably equivalent to $Y$ given by 
\[
Y(V) = \prod_{\text{$n$-cells of $B$}}\left\{ s:G/H\times D^n \to X(V)\mid \pi_{X(V)}\circ s = 1 \text{ and } s|_{G/H\times \partial D^n} = i_{X(V)} \right\}
\]
and bonding maps given by applying the bonding maps of $Fj$ to sections restricted to a single $n$-cell. The equivalence $Y\to Fj$ is given by sending a section on $G/H\times D^n\to X(V)$ to a section on $B^{(n)}\to X(V)$ which agrees with the original section on $G/H\times D^n$ and is the zero section everywhere else.

Now observe that since $G/H\times D^n$ is equivariantly contractible, sections $G/H\times D^n\to X(V)$ can, without loss of generality, be taken to be sections $G/H\times D^n\to (G/H\times D^n)\times X(V)_b$ where $X(V)_b$ is the fiber over $b$. Since we have assumed $X$ to be fibrant (or else fibrantly replaced it), there is a canonical choice of fiber spectrum of $X$. Observe that the space of sections of $G/H\times D^n\times X(V)_b$ over $G/H\times D^n$ is homeomorphic to the space of maps from $G/H\times D^n\to X(V)_b$. So we may take $Y(V)$ to be
\[
Y(V)\simeq \prod_{\text{$n$-cells of $B$}}\left\{ s:G/H\times D^n \to X(V)_b\mid s|_{G/H\times \partial D^n} = i_{X(V)} \right\}.
\]
This is equivalent to 
\[
\prod_{\text{$n$-cells of $B$}}\left\{ s:G/H\wedge_+ S^n \to X(V)_b\right\}.
\]
The same argument shows that $Fl$ is equivalent to $Z$ given by 
\[
Z(V) = \prod_{\text{$n$-cells of $B$}} \left\{ G/H\wedge_+ S^n\to (A\barsmash X)_b \right\}.
\]
Now since products and wedges of spectra are stably equivalent, and smashing with $A$ commutes with wedge sums of spectra, we can see that 
\[
A\wedge Y\simeq Z.
\]
Moreover, after passing through the identifications
\[
A\wedge Y \simeq A\wedge Fj \qquad\text{and}\qquad Z\simeq Fl
\]
this stable equivalence becomes precisely the map $A\wedge Fj\to Fl$ in the diagram above which is a stable equivalence. 

Since we have a map of fiber sequences of spectra where the left and right vertical maps are equivalences, the long exact sequence of homotopy groups and the five lemma show that the middle map 
\[
A \wedge\Gamma_{B^{(n)}}(X)\to \Gamma_{B^{(n)}}(A\barsmash X)
\]
is also a stable equivalence. We conclude that
\[
A\wedge \Gamma_B(X)\simeq \Gamma_B(A\barsmash X).
\]
\end{proof}

\begin{lem}\label{lem:fibrant concentrated isotropy}
Let $X\in G\Osp(B)$. If $X$ is stably fibrant, level $h$-cofibrant, level $h$-fibrant, and has isotropy concentrated in $H$, then $\widetilde{EP_H}\barsmash X$ is fibrant as well.
\end{lem}

\begin{proof}
We begin by noting that $X$ is fibrant if it is level $q$-fibrant and its underlying orthogonal spectrum is an $\Omega$-spectrum. We may choose $\widetilde{EP_H}$ to be $h$-cofibrant since it was only defined up to homotopy. Then since $X$ is level $h$-cofibrant and level $h$-fibrant, we may conclude that $\widetilde{EP_H}\barsmash X$ is level $h$-fibrant as well. It is therefore automatically level $q$-fibrant. 

It remains to show that $\widetilde{EP_H}\barsmash X$ is an $\Omega$-spectrum. Since $X$ is an $\Omega$-spectrum, the bonding maps $X(V)\to \Omega^{W-V}X(W)$ are weak equivalences. The bonding maps in $\widetilde{EP_H}\barsmash X$ are $\widetilde{EP_H}\barsmash X(V) \to \Omega^{W-V}(\widetilde{EP_H}\barsmash X(W))$. We recover the original bonding maps on the $H$-fixed points, whereas everywhere else both spaces are contractible. Therefore, our new bonding maps are weak equivalences as desired. 
\end{proof}

\begin{lem}\label{lem:pullbacks commute with fixed points}
Let $\phi:A\to B$. The right derived pullback functor $\phi^*:G\Osp(B)\to G\Osp(A)$ commutes with the left derived geometric fixed points.
\end{lem}

\begin{proof}
It suffices to prove the result on the full subcategory on the objects in $G\Osp(B)$ which are stably fibrant, level $h$-fibrant and level $h$-cofibrant.

Using \cref{lem:fibrant concentrated isotropy} and the fact that $\Phi^H$ is monoidal, we see that
\[
\Phi^H(X)\simeq (\widetilde{EP_H}\barsmash f(\widetilde{EP_H}\barsmash X))^H
\]
where $f(-)$ is the stable fibrant replacement. Moreover, this is a weak equivalence between fibrant objects so we obtain
\[
\phi^*(\Phi^H(X)) \simeq \phi^*((\widetilde{EP_H}\barsmash f(\widetilde{EP_H}\barsmash X))^H).
\]
By \cite[Lemma 7.1.2]{malkiewich_convenient} and \cite[Proposition 4.5.1]{malkiewich_convenient}, $\phi^*$ commutes with categorical fixed points and $\widetilde{EP_H}\barsmash(-)$ respectively. So that the above equation becomes 
\begin{align*}
\phi^*(\Phi^H(X)) 
& \simeq \phi^*((\widetilde{EP_H}\barsmash f(\widetilde{EP_H}\barsmash X))^H)\\
& \simeq (\widetilde{EP_H}\barsmash \phi^*f(\widetilde{EP_H}\barsmash X))^H\\
& \simeq \Phi^H(\phi^*(f(\widetilde{EP_H}\barsmash X)))
\end{align*}
where the last line again follows from \cref{lem:fibrant concentrated isotropy}. We now claim that the canonical map 
\[
\Phi^H(\phi^*(\widetilde{EP_H}\barsmash X))\rightarrow \Phi^H(\phi^*(f(\widetilde{EP_H}\barsmash X)))
\]
is a stable equivalence. Observe that $\widetilde{EP_H}\barsmash X\to f(\widetilde{EP_H}\barsmash X)$ is an acyclic cofibration. Since $\phi^*$ is left Quillen it preserves acyclic cofibrations. By \cite{mp2}, $\Phi^H(-)$ also preserves acyclic cofibrations. So the above map is a weak equivalence as desired. Since $\phi^*$ commutes with $\widetilde{EP_H}\barsmash (-)$ and $\Phi^H(-)$ is monoidal,
\[
\Phi^H(\phi^*(\widetilde{EP_H}\barsmash X))\simeq \Phi^H(\phi^*(X)).
\]
Finally, combining the above equations yields
\[
\phi^*(\Phi^H(X)) \simeq \Phi^H(\phi^*(X))
\]
as desired. 
\end{proof}

\begin{prop}\label{prop:sections commute with geom fp}
Let $B$ be a finite $G$-CW-complex. Let $X\in G\Osp(B)$ be a $G$-spectrum parameterized over $B$ and $H$ a closed subgroup of $G$. The right derived spectrum of sections, $\R\Gamma_B(-)$, commutes with the left derived geometric fixed points, $\Phi^H(-)$. 
\end{prop}

\begin{proof}
It suffices to show the result on the full subcategory of $G\Osp(B)$ on the objects which are stably fibrant, level $h$-fibrant and level $h$-cofibrant.

We proceed by induction on the skeleta of $B$. We first show the base case. Suppose that $B$ is some finite $G$-set. We wish to show that
\[
\Gamma_{B^H}\Phi^HX  \simeq \Phi^H \Gamma_B X.
\]
We first observe that non-equivariantly
\[
\Gamma_BX = \prod_{b\in B} X_b
\]
To make sense of this statement equivariantly, we define a $G$ action on $\prod_{b\in B} X_b$ by 
\[
g\cdot(x_1,\dots,x_k)  = (gx_{g^{-1}1},\dots,gx_{g^{-1}k})
\]
where 
\[
(x_1,\dots,x_k) \in X_1(V)\times\dots\times X_k(V)
\]
for some $G$-representation $V$. 

This action makes the unit and counit maps from the non-equivariant adjunction between pullback to $B$ and sections on $B$ into equivariant maps defining an adjunction in   $G\Osp(B)$. 

Now observe that  
\begin{align*}
\Phi^H\left(\prod_{b\in B} X_b\right) \simeq \prod_{b\in B^H} \Phi^H(X_b).
\end{align*}

Moreover, if $b\in B^H$ then the fiber spectrum over $b$ is the pullback of $X$ along $\{b\}\to B$. Thus, $\Phi^H(X_b)\cong (\Phi^H X)_b$ by \cref{lem:pullbacks commute with fixed points}. This proves the claim for sections on a finite $G$-set.

We now proceed to the inductive step. Let $B$ be some finite $G$-CW complex and suppose that 
\[
\Gamma_{B^{(n-1)^H}} \Phi^H X \simeq \Phi^H \Gamma_{B^{(n-1)}} X
\]
Again by \cref{lem:serre fibration}, we have a levelwise Serre fibration
\[
\Gamma_{B^{(n)}} X\xrightarrow{j}\Gamma_{B^{(n-1)}}X 
\]
The functor $\Phi^H$ is excisive \cite{alaska} and thus preserves fiber sequences. In other words, we have a fiber sequence
\[
\Phi^H Fj\to \Phi^H\Gamma_{B^{(n)}}X\xrightarrow{\Phi^H (j)}\Phi^H \Gamma_{B^{(n-1)}}X.
\]
We also have a fiber sequence
\[
Fl\to\Gamma_{B^{(n)^H}}\Phi^HX \xrightarrow{l}\Gamma_{B^{(n-1)^H}}\Phi^HX
\]
where $l$ is a levelwise Serre fibration by \cref{lem:serre fibration}. 

We now build a map between these two fiber sequences. By \cref{lem:pullbacks commute with fixed points}, there is a natural isomorphism commuting the pullback functor with geometric fixed points. Composing this with the unit and the counit of the adjunction between pullback and $ \Gamma_B$ we obtain a natural transformation from $\Phi^H\Gamma_B$ to $\Gamma_B\Phi^H$. These maps assemble into a commuting square:
\begin{figure}[H]
\center
\begin{tikzcd}
\Phi^H\Gamma_{B^{(n)}} X \arrow[r,"\Phi^H j"]\arrow[d]& \Phi^H\Gamma_{B^{(n-1)}} X\arrow[d]\\
\Gamma_{B^{(n)^H}} \Phi^H X \arrow[r,"l"]& \Gamma_{B^{(n-1)^H}}\Phi^H X.
\end{tikzcd}
\end{figure}

Then the map induced on the fibers turns this into a map of fiber sequences:
\begin{diagram}\label{dia:fiber sequence}
\Phi^H Fj\arrow[d]\arrow[r] &\Phi^H\Gamma_{B^{(n)}} X \arrow[r,"\Phi^H j"]\arrow[d]& \Phi^H\Gamma_{B^{(n-1)}} X\arrow[d]\\
Fl\arrow[r]&\Gamma_{B^{(n)^H}} \Phi^H X \arrow[r,"l"]& \Gamma_{B^{(n-1)^H}}\Phi^H X. 
\end{diagram}

We now argue that the left vertical map is a stable equivalence. As in \cref{prop:sections commutes with smash}, we may rewrite $Fj$ as
\[
Fj\simeq \prod_{\text{$n$-cells of $B$}}F(G/K\wedge_+ S^n ,X_b)
\]
where $F(-,-)$ is the cotensor spectrum. Now, since $\Phi^H(-)$ commutes with wedge sums, it commutes with finite products. Additionally, since $G/K\wedge_+ S^n$ is a finite CW-complex (and thus dualizable), the canonical map $\Phi^H:F(G/K\wedge_+ S^n,X)\to F((G/K)^H\wedge_+ S^n,\Phi^H(X)_b)$ is a weak equivalence. Finally, since $\Phi^H(-)$ preserves fiber sequences, and thus cofiber sequences, we observe that $\Phi^H(X)_b\simeq \Phi^H(X_b)$. Putting these three facts together, we see that

\begin{align*}
\Phi^H(Fj)
& \simeq \Phi^H\left(\prod_{\text{$n$-cells of $B$}}F(G/K\wedge_+ S^n ,X_b)\right)\\
& \simeq \prod_{\text{$n$-cells of $B^H$}}F((G/K)^H\wedge_+ S^n,\Phi^H(X)_b).
\end{align*} 

Similarly, $Fl$ is stably equivalent to  
\[
\prod_{\text{$n$-cells of $B^H$}}F((G/K)^H\wedge_+ S^n\to \Phi^H (X)_b.
\]

Moreover, the map $\Phi^H(Fj)\to Fl$ of \cref{dia:fiber sequence} is precisely this chain of weak equivalences.

Finally, the long exact sequence of homotopy groups associated to a fiber sequence and the five-lemma show that
\[
\Gamma_{B^{(n)^H}} \Phi^H X\simeq \Phi^H \Gamma_{B^{(n)}} X
\]
as desired. 
\end{proof}

\bibliographystyle{amsalpha}
\bibliography{streamlined}

\end{document}

%% file: preamble_defs.tex
\newcommand{\sJ}{\mathscr{J}}

\newcommand{\cS}{\mathcal{S}}

\newcommand{\R}{\mathbb R}

\newcommand{\Sph}{\mathbb{S}}

\newcommand{\cat}[1]{\textup{\textbf{{#1}}}}

\newcommand{\GTop}{\cat{GTop}}

\newcommand{\All}{\textup{All}}
\newcommand{\free}{\textup{free}}
\newcommand{\id}{\textup{id}}

\newcommand{\colim}{\textup{colim}\,}
\newcommand{\hocolim}{\textup{hocolim}\,}

\newcommand{\barsmash}{\,\overline\wedge\,}

\newcommand{\tr}{\textup{tr}\,}

\DeclareMathOperator{\res}{\text{res}}

\newcommand{\Osp}{\mathcal{OS}}

\newcommand{\po}[2]{\ar@{}@<{#2}>[rd]|({#1})*\txt{\Large $\ulcorner$}}
\newcommand{\pb}[2]{\ar@{}@<{#2}>[rd]|({#1})*\txt{\Large $\lrcorner$}}

\DeclareMathOperator{\Lie}{Lie}

%% file: preamble_refs_angelica.tex
\usepackage[unicode=true, pdfusetitle,
 bookmarks=true,bookmarksnumbered=false,
 breaklinks=false,
 backref=false,
 colorlinks=true,
 linkcolor=blue,
 citecolor=blue,
 urlcolor=blue,
 final
]{hyperref}
\usepackage{amsthm}
\usepackage{aliascnt}
\usepackage{import}
\usepackage[capitalise]{cleveref}

\newcommand{\newrefformat}[2]{}

\newenvironment{diagram}{\crefalias{equation}{diagram}\begin{equation}\begin{tikzcd}}{\end{tikzcd}\end{equation}}
\newenvironment{diagram*}{\crefalias{equation*}{diagram}\begin{equation*}\begin{tikzcd}}{\end{tikzcd}\end{equation*}}

\crefname{diagram}{Diagram}{Diagrams}
\Crefname{diagram}{Diagram}{Diagrams}

\theoremstyle{plain}   
\newtheorem{thm}{Theorem}[section] 
\makeatletter\let\c@thm\c@thm\makeatother
\newtheorem{cor}[thm]{Corollary}
\makeatletter\let\c@cor\c@thm\makeatother
\newtheorem{lem}[thm]{Lemma}
\makeatletter\let\c@lemma\c@thm\makeatother
\newtheorem{prop}[thm]{Proposition}
\makeatletter\let\c@prop\c@thm\makeatother

\makeatletter\let\c@claim\c@thm\makeatother

\theoremstyle{definition}

\newtheorem{df}[thm]{Definition}
\makeatletter\let\c@defn\c@thm\makeatother

\makeatletter\let\c@const\c@thm\makeatother
\newtheorem{notn}[thm]{Notation}
\makeatletter\let\c@notn\c@thm\makeatother

\makeatletter\let\c@outline\c@thm\makeatother

\makeatletter\let\c@propty\c@thm\makeatother

\makeatletter\let\c@problem\c@thm\makeatother

\makeatletter\let\c@conj\c@thm\makeatother

\theoremstyle{remark}

\newtheorem{rmk}[thm]{Remark}
\makeatletter\let\c@rem\c@thm\makeatother
\newtheorem{ex}[thm]{Example}
\makeatletter\let\c@ex\c@thm\makeatother

\makeatletter\let\c@observationn\c@thm\makeatother

\makeatletter
\let\c@equation\c@thm
\numberwithin{equation}{section}
\makeatother

\crefname{lemma}{Lemma}{Lemmas}
\crefname{thm}{Theorem}{Theorems}
\crefname{defn}{Definition}{Definitions}
\crefname{notn}{Notation}{Notations}
\crefname{const}{Construction}{Constructions}
\crefname{prop}{Proposition}{Propositions}
\crefname{rem}{Remark}{Remarks}
\crefname{cor}{Corollary}{Corollaries}
\crefname{equation}{Equation}{Equations}
\crefname{ex}{Example}{Examples}
\crefname{propty}{Property}{Properties}
\crefname{problem}{Problem}{Problems}

\hypersetup{linktoc=all}

%% file: streamlined.bib
@book{dhks,
AUTHOR = {Dwyer, William G. and Hirschhorn, Philip S. and Kan, Daniel M.
and Smith, Jeffrey H.},
TITLE = {Homotopy limit functors on model categories and homotopical
categories},
SERIES = {Mathematical Surveys and Monographs},
VOLUME = {113},
PUBLISHER = {American Mathematical Society, Providence, RI},
YEAR = {2004},
PAGES = {viii+181},
ISBN = {0-8218-3703-6},
MRCLASS = {18G55 (18A30 18B99 55U35)},
MRNUMBER = {2102294},
MRREVIEWER = {Beno\^{\i}t Fresse},
}

@article{elmendorf83,
title={Systems of fixed point sets},
author={Elmendorf, T.},
journal={Trans. Amer. Math. Soc.},
fjournal={Transaction of the American Mathematical Society},
volume={277},
year={1983},
number={1},
pages={275-284},
url={https://www.jstor.org/stable/1999356},
}

@article{gimm,
	author = {Goodwillie, T.~G. and Igusa, K. and Malkiewich, M. and Merling, M.},
	title = {On the functoriality of the space of equivariant smooth $h$-cobordisms},
	journal = {submitted},
	year = {2023},
	url={https://doi.org/10.48550/arXiv.2303.14892},
}

@article{Hausschild1974,
	author = {Hausschild, Henning},
	journal = {Mathematische Zeitschrift},
	language = {ger},
	pages = {165-172},
	title = {Bordismentheorie stabil gerahmter {$G$}-Mannigfaltigkeiten.},
	url = {http://eudml.org/doc/172119},
	volume = {139},
	year = {1974},
}

@article{heath_kamps,
	author={Heath, P. and Kamps, K.},
	title={Note on attaching {D}old fibrations},
	journal={Canadian Mathematical Bulletin},
	volume={21},
	issue={3},
	year={1978},
	pages={365-367}
}

@book{hhr_book, 
place={Cambridge},
 series={New Mathematical Monographs}, 
title={Equivariant Stable Homotopy Theory and the {K}ervaire Invariant Problem}, 
publisher={Cambridge University Press}, 
author={Hill, Michael A. and Hopkins, Michael J. and Ravenel, Douglas C.}, 
year={2021},
 collection={New Mathematical Monographs},
}

@article {hhr,
	AUTHOR = {Hill, M. A. and Hopkins, M. J. and Ravenel, D. C.},
	TITLE = {On the nonexistence of elements of {K}ervaire invariant one},
	JOURNAL = {Ann. of Math. (2)},
	FJOURNAL = {Annals of Mathematics. Second Series},
	VOLUME = {184},
	YEAR = {2016},
	NUMBER = {1},
	PAGES = {1--262},
	ISSN = {0003-486X},
	MRCLASS = {55P91 (55N22 55P42 55Q45 55T15 55U35 57R15)},
	MRNUMBER = {3505179},
	MRREVIEWER = {Paul G. Goerss},
	URL = {https://doi.org/10.4007/annals.2016.184.1.1},
}

@article{kervaire-milnor,
 ISSN = {0003486X, 19398980},
 URL = {http://www.jstor.org/stable/1970128},
 author = {Michel A. Kervaire and John W. Milnor},
 journal = {Annals of Mathematics},
 number = {3},
 pages = {504--537},
 publisher = {[Annals of Mathematics, Trustees of Princeton University on Behalf of the Annals of Mathematics, Mathematics Department, Princeton University]},
 title = {Groups of Homotopy Spheres: I},
 urldate = {2024-11-21},
 volume = {77},
 year = {1963}
}

@article{kosniowski_bordism,
	title={Equivariant Stable Homotopy and Framed Bordism},
	author={Kosniowski, Czes},
	journal={Transactions of the American Mathematical Society},
	volume={219},
	pages={225--234},
	year={1976},
	publisher={American Mathematical Society}
}

@incollection {luck_surgery,
	AUTHOR = {L{\"u}ck, Wolfgang},
	TITLE = {A basic introduction to surgery theory},
	BOOKTITLE = {Topology of high-dimensional manifolds, {N}o. 1, 2 ({T}rieste,
	2001)},
	SERIES = {ICTP Lect. Notes},
	VOLUME = {9},
	PAGES = {1--224},
	PUBLISHER = {Abdus Salam Int. Cent. Theoret. Phys., Trieste},
	YEAR = {2002},
	MRCLASS = {57R67 (57R65)},
	MRNUMBER = {1937016},
}

@article{malkiewich_convenient,
	title = {A convenient category of parameterized spectra},
	author = {Malkiewich, Cary},
	journal = {arXiv preprint https://arxiv.org/abs/2305.15327},
	year = {2023}	
}

@article{mp2,
    AUTHOR = {Malkiewich, Cary and Ponto, Kate},
     TITLE = {Periodic points and topological restriction homology},
   JOURNAL = {Int. Math. Res. Not. IMRN},
  FJOURNAL = {International Mathematics Research Notices. IMRN},
      YEAR = {2022},
    NUMBER = {4},
     PAGES = {2401--2459},
      ISSN = {1073-7928},
   MRCLASS = {57Q10},
  MRNUMBER = {4381922},
       DOI = {10.1093/imrn/rnaa174},
       URL = {https://doi.org/10.1093/imrn/rnaa174},
}

@article{mandell2002equivariant,
    AUTHOR = {Mandell, M. A. and May, J. P.},
    TITLE = {Equivariant orthogonal spectra and {$S$}-modules},
    JOURNAL = {Mem. Amer. Math. Soc.},
    FJOURNAL = {Memoirs of the American Mathematical Society},
    VOLUME = {159},
    YEAR = {2002},
    NUMBER = {755},
    PAGES = {x+108},
    ISSN = {0065-9266},
    MRCLASS = {55P91 (18E30 55P42 55P43 55P48)},
    MRNUMBER = {1922205},
    MRREVIEWER = {J. P. C. Greenlees},
    DOI = {10.1090/memo/0755},
    URL = {https://doi.org/10.1090/memo/0755},
}

@book{alaska,
	AUTHOR = {May, J.P.},
	TITLE = {Equivariant Homotopy and Cohomology Theory},
	SERIES = {Regional Conference Series in Mathematics},
	PUBLISHER = {American Mathematical Society},
	YEAR = {1996},
	ISBN = {0821803190},
}

@book{ms,
    AUTHOR = {May, J. P. and Sigurdsson, J.},
    TITLE = {Parametrized homotopy theory},
    SERIES = {Mathematical Surveys and Monographs},
    VOLUME = {132},
    PUBLISHER = {American Mathematical Society, Providence, RI},
    YEAR = {2006},
    PAGES = {x+441},
    ISBN = {978-0-8218-3922-5; 0-8218-3922-5},
    MRCLASS = {55P42 (19L99 55N20 55N22 55P91)},
    MRNUMBER = {2271789},
    MRREVIEWER = {A. A. Ranicki},
    DOI = {10.1090/surv/132},
    URL = {https://doi.org/10.1090/surv/132},
}

@book{milnor_diff_top,
	AUTHOR = {Milnor, John},
    TITLE = {Topology From the Differentiable Viewpoint},
    SERIES = {Princeton Landmarks in Mathematics},
    PUBLISHER = {Princeton University Press, Princeton, New Jersey},
    YEAR = {1965},
    ISBN = {0813901812},
}

@incollection {segal1970equivariant,
    AUTHOR = {Segal, G. B.},
     TITLE = {Equivariant stable homotopy theory},
 BOOKTITLE = {Actes du {C}ongr\`es {I}nternational des {M}ath\'ematiciens
              ({N}ice, 1970), {T}ome 2},
     PAGES = {59--63},
 PUBLISHER = {Gauthier-Villars \'Editeur, Paris},
      YEAR = {1971},
   MRCLASS = {55D40 (55E10 57D85 57E15)},
  MRNUMBER = {423340},
MRREVIEWER = {Henning\ Hauschild},
}

@book {schwede_global,
	AUTHOR = {Schwede, Stefan},
	TITLE = {Global homotopy theory},
	SERIES = {New Mathematical Monographs},
	VOLUME = {34},
	PUBLISHER = {Cambridge University Press, Cambridge},
	YEAR = {2018},
	PAGES = {xviii+828},
	ISBN = {978-1-108-42581-0},
	MRCLASS = {55P42 (18G55 19D99 55P91 55U35)},
	MRNUMBER = {3838307},
	MRREVIEWER = {Gregory Z. Arone},
	DOI = {10.1017/9781108349161},
	URL = {https://doi.org/10.1017/9781108349161},
}

@book{stong68,
	title={Notes on Cobordism Theory},
	author = {Stong, Robert},
	year={1968},
	publisher={Princeton University Press}

}

@article{tomDieckEquivariantBordism,
author={tom dieck, Tammo},
title={Orbittypen und aquivariante Homologie II.},
journal={Arch. Math.},
volume={26},
year={1975},
number={6},
pages={650-662},
}

@article{manifold_approach,
  title={Algebraic {$K$}-theory of spaces, a manifold approach},
  author={Waldhausen, Friedhelm},
  journal={Current trends in algebraic topology, Part},
  volume={1},
  pages={141--184},
  year={1982},
}

@article {waner_1984_eq_bordism,
	author = {Waner, Stefan},
	title = {Equivariant {$RO(G)$}-graded Bordism Theories},
	journal = {Topology and its Applications},
	volume = {17},
	issue={1},
	year = {1984},
	pages = {1--26},
	url = {https://doi.org/10.1016/0166-8641(84)90019-1}

}

@article{wang_xu_61stem,
	author={Wang, Guozhen and Xu, Zhouli},
	title={The Triviality of the 61-stem in the Stable Homotopy Group of Spheres},
	journal={Annals of Mathematics},
	volume={186},
	issue={2},
	year={2017},
	pages={501--580},
	url={https://doi.org/10.4007/annals.2017.186.2.3}	
}

@article{wasserman,
title = {Equivariant differential topology},
journal = {Topology},
volume = {8},
number = {2},
pages = {127-150},
year = {1969},
issn = {0040-9383},
doi = {https://doi.org/10.1016/0040-9383(69)90005-6},
url = {https://www.sciencedirect.com/science/article/pii/0040938369900056},
author = {Arthur G. Wasserman}
}
